\documentclass[11pt]{article}
\usepackage[margin=1in]{geometry}
\usepackage{amsmath, amssymb, graphicx, float, amsthm}

\newtheorem{lemma}{Lemma}

\theoremstyle{definition}

\newtheorem{remark}{Remark}

\usepackage{subcaption}

\theoremstyle{remark}

\newtheorem{proposition}{Proposition}

\usepackage{algorithm}
\usepackage[noend]{algpseudocode}
\usepackage{multirow}
\usepackage{tabularx}
\usepackage{makecell}
\usepackage{booktabs}
\usepackage{natbib}
\bibpunct{(}{)}{;}{a}{}{,}
\usepackage{tikz}
\usetikzlibrary{arrows.meta, positioning, fit, shapes.arrows}
\usepackage[figuresright]{rotating}
\usepackage{url}
\usepackage{authblk}
\usepackage[colorlinks=true, urlcolor=blue, citecolor=black, linkcolor=black]{hyperref}
\begin{document}

\title{The operational impact of registry size, cycle length and blood type distribution in Multi-Registry Kidney Exchange Programs}

\author[a,*]{Utkarsh Verma}
\author[b]{Narayan Rangaraj}
\author[c]{Viswanath Billa}
\author[d]{Deepa Usulumarty}

\affil[a]{Department of Operations Management and Information System, Schulich School of Business, York University, Toronto M3J 1P3, Canada}
\affil[b]{Department of Industrial Engineering and Operations Research, Indian Institute of Technology Bombay, Mumbai, Maharashtra 400076, India}
\affil[c]{Department of Nephrology, Apex Kidney Care, Sushrut Hospital, Mumbai, Maharashtra 400071, India}
\affil[d]{Department of Nephrology, Bombay Hospital, Mumbai, Maharashtra 400020, India}

\date{}
\maketitle

\begin{center}
\small
\textsuperscript{*}Corresponding author: Utkarsh Verma (\href{mailto:utkarshv@schulich.yorku.ca}{utkarshv@schulich.yorku.ca}), ORCID: \href{https://orcid.org/0000-0003-4756-6856}{0000-0003-4756-6856}
\end{center}

\begin{abstract}
Kidney exchange programs (KEPs) address donor--recipient incompatibility by exchanging donors
between incompatible pairs, but single-center KEPs often suffer from limited pool sizes that
reduce matching efficiency. Multi-registry kidney exchange programs (mKEPs) offer a promising
solution but face challenges including heterogeneous constraints across registries, cycle-length
bounds, and data-sharing limitations. This study uses simulation to compare mKEP allocation
against individual registry operation, contrasting unconstrained pooling with two safeguarded
mechanisms: a cumulative individual-rationality guarantee and a Shapley-value-based fair-share
mechanism. We examine how registry size, blood-type distribution, cycle-length bounds, and
dropout probability affect the size and distribution of the achievable gain. Our central finding
is that the benefit of joining an mKEP is systematically uneven: under blood-group-composition
asymmetry, the easier-to-match registry gains less than its partner even under both safeguarded
mechanisms, and a larger registry gains less than a smaller one when pooled. A registry combining
a larger arrival rate with an easier-to-match pool may see its gain under either safeguarded
mechanism fall too small to be practically significant, while registries with a higher dropout
rate see a larger benefit from pooling. This
unevenness concentrates in O-type recipients: an easier-to-match registry sees fewer O-type
transplants within itself under pooling, though O-type transplant rates rise system-wide. Tighter
cycle-length bounds increase, rather than diminish, the relative transplant-volume benefit of
pooling, while match quality is only modestly affected by any factor examined. These results
highlight the importance of pairing multi-registry collaboration with a carefully designed,
equitable benefit-sharing mechanism to keep participation attractive for all registries involved.
\end{abstract}



\section{Introduction}\label{intro}

Over the past two decades, Kidney Exchange Programs (KEPs) have evolved into a vital mechanism for increasing transplant opportunities by enabling incompatible donor--recipient pairs to exchange donors in order to achieve compatibility \citep{Roth2004,Roth2005,Abraham2007,Wallis2011,Constantino2013,Anderson2014,Ferrari2015,Kute2018, Ashlagi2024, Barkel2025}. KEPs initially focused exclusively on simple swap exchanges, where two incompatible donor--recipient pairs exchange donors to achieve mutual compatibility \citep{Roth2004}. Over time, these programs evolved to include more complex cycle exchanges, involving three or more pairs participating in closed exchange cycles, thereby increasing the flexibility and match potential within the system.

The integration of non-directed donors including altruistic and deceased donors, marked a significant advancement in KEPs \citep{Verma2023}. These donors can initiate exchange chains, sequences of transplants that substantially expand the pool of transplants and enhance the efficiency of kidney allocation \citep{Veale2009,Melcher2013}. There is inclusion of ABO-incompatible transplants within KEPs been proposed as a means of further widening match opportunities \citep{Karami2019}. In parallel, it has become increasingly common for compatible pairs to voluntarily join exchange registries in pursuit of better immunological or demographic matches, particularly to improve the recipient's expected transplant outcome \citep{Weng2017}. Recent advances in immunology have demonstrated the potential for enzymatic conversion of donor organs to blood type O, thereby mitigating ABO incompatibility in kidney paired donation (KPD) \citep{Yang2026}.

In recent years, kidney exchange programs have experienced significant growth, particularly through the development of international and global exchange models \citep{Biro2019, Sharifi2025}. In Europe, international collaborations have become increasingly common, beginning with the first cross-border kidney exchange transplants between the Czech Republic and Austria in 2016 \citep{Biro2021}. Since then, two additional international programs have been actively operating: the Spain--Italy--Portugal partnership, and the Scandinavian Transplant Exchange Program (STEP), which involves Sweden, Norway, and Denmark under the Scandiatransplant framework \citep{Biro2017,Biro2019}.

Beyond regional collaborations, the concept of Global Kidney Exchange (GKE) has emerged, aiming to facilitate transplants between donor--recipient pairs from low- or middle-income countries and those from high-income countries \citep{Rees2022}. Under this model, the transplant costs for pairs from lower-income countries are covered by their high-income partners, thereby enabling broader access to compatible matches \citep{Pullen2017,Rees2017}. While a limited number of such transplants have been successfully conducted \citep{Kute2017}, the ethical implications of GKE remain a topic of ongoing debate, with concerns raised regarding equity, exploitation, and long-term sustainability \citep{Ambagtsheer2020}.

In India, both single-center kidney exchange registries and multi-center collaborations involving combined hospital-based registries are currently operational \citep{Kute2017, Billa2018}. However, these efforts remain limited in their ability to significantly alleviate the burden of kidney transplantation, primarily due to small registry sizes and logistical constraints. As such, the development of multi-registry or multi-hospital kidney exchange programs represents a promising next step toward enhancing the efficiency and reach of exchange-based transplant solutions.

Given the diversity in state-level policies governing kidney exchange transplants across India, inter-state collaborations may face challenges similar to those encountered in international exchanges across European countries. This makes the Indian context particularly suitable for evaluating mechanisms developed for cross-state exchanges. Furthermore, multi-registry exchange programs in India could also incorporate networks of private hospitals, further expanding the pool of donor--recipient pairs and increasing the potential for high-quality matches.

These multi-registry kidney exchange programs (mKEPs) face several practical and strategic challenges, including issues related to data sharing, varying constraints on cycle lengths, and decisions about whether transplants should be executed simultaneously or non-simultaneously \citep{Mincu2020, Colley2025}. A critical concern in such collaborations is the potential for individual registries to have strategic incentives to withhold complete information. This reluctance to share full data can hinder the natural formation of coalitions, as demonstrated by \citet{Hajaj2015}, who show that self-interest can prevent registries from participating fully in joint allocation systems.

Two primary models for data sharing in mKEPs have been implemented internationally: sequential merger and full merger. In the sequential merger model, registries first attempt to match pairs internally and then share the unmatched pairs with partner registries. A prominent example of this approach is the Spain--Italy--Portugal collaboration, where each country completes internal matches before engaging in cross-border matching \citep{Mincu2020}. In contrast, the full merger model involves complete data sharing across registries, enabling algorithms to optimize matches across the entire joint pool. This model has been adopted in the collaboration between the Czech Republic and Austria, as well as in the United States, where multi-hospital kidney exchange programs were initiated in 2008, though certain conditions and restrictions still apply for specific patient categories \citep{Hanto2008,Viklicky2020}. When a match fails, the literature has proposed various recourse strategies, notably internal recourse and subset recourse, with evidence suggesting that subset recourse yields superior outcomes compared to internal recourse \citep{Smeulders2022, Matyasi2024}.

Both sequential and full merger approaches to mKEPs present inherent limitations. In a sequential merger, suboptimality can arise because certain donor--recipient pairs are matched and removed before data sharing occurs with other registries, thereby reducing the overall efficiency of the joint system. Conversely, in a full merger -- where all registries share complete information -- the resulting allocation may disproportionately benefit some registries over others, raising concerns around fairness and, in the extreme, around whether a registry is left worse off than it would have been alone.

Two distinct lines of work have addressed this concern. The first anchors the round's allocation to an individual-rationality guarantee directly: \citet{Hajaj2015} propose a dynamic credit mechanism that is simultaneously strategy-proof and, in the long run, individually rational, each transplant center is guaranteed at least as many matches as it could have secured acting alone; \citet{Ashlagi2011,Ashlagi2012} study individual rationality as a hard participation constraint in large-scale multi-hospital kidney exchange, characterizing the efficiency price of guaranteeing that no hospital is made worse off by pooling; and \citet{Roth2014} show that free-riding incentives can still undermine participation even when a nominal individual-rationality guarantee is in place. The second line anchors the round's allocation to an axiomatic, cooperative-game fair-share concept instead: \citet{Klimentova2021} propose a lexicographic, multi-round credit-balancing mechanism comparing several solution concepts the benefit value, the contribution value, the Shapley value, and the nucleolus, showing that balancing consecutive matching rounds, rather than solving each in isolation, produces markedly more equitable long-run outcomes without depending on the arbitrary order in which registries' data is processed; \citet{Biro2020} compare a benefit-value-based target against a Shapley-value-based one directly and find the latter yields a lower standard deviation in realized benefit across registries; and \citet{Benedek2024} study a closely related credit-based, deviation-minimizing mechanism at considerably larger scale up to fifteen countries and confirming, independently and at that scale, that a Shapley-value-based target outperforms several alternative cooperative-game solution concepts.

These two lines of work, however, have not, to our knowledge, been directly compared within a single simulation framework: it remains unclear whether an individual-rationality-anchored mechanism and an axiomatic fairness-anchored mechanism produce materially different outcomes for a given registry, nor whether a mechanism's decision to permanently retain a registry's surplus (as in an individual-rationality guarantee that is never revised downward) versus redistributing it in later rounds (as in a signed, budget-balancing credit) matters in practice. Separately, while the effect of genuinely heterogeneous registry or country parameters on multi-agent KEP outcomes has been examined using real, country-specific arrival volumes and blood-type distributions \citep{Druzsin2024}, that work is directed at comparing collaboration policies and optimization criteria across heterogeneous pools, rather than at how such heterogeneity interacts with the \emph{choice} of fairness or individual-rationality mechanism, in particular, whether a registry with an easier-to-match patient pool retains its surplus indefinitely under one mechanism but has it redistributed under another, and how large registry heterogeneity must become before the achievable gain from joining an mKEP is too small, or too asymmetrically distributed, to justify a registry's continued participation at all.

Key questions therefore include whether registries of differing sizes, blood-group compositions, or operational risk profiles would be willing to collaborate in a joint allocation system when the distribution of benefits is potentially disproportionate, and whether a simple individual-rationality is sufficient to sustain that participation or whether a more demanding fairness criterion and specifically, a criterion that redistributes rather than permanently cedes surplus is required. For example, in contexts like India, blood group distributions and registry behaviors vary significantly: some registries may enroll O-type donors and AB-type recipients predominantly in cases of cross-match failure or low compatibility scores, while others may not frequently enroll compatible pairs. This heterogeneity, together with India's diversity of state-level policy and registry scale, remains largely unexamined in a literature that is otherwise calibrated overwhelmingly to European and North American registries \citep{Druzsin2024,Ashlagi2012}.

Additionally, differences in operational constraints such as allowable cycle lengths or dropout probabilities raise concerns about the willingness of more stable or resource-rich registries to participate in joint systems with higher-risk or more constrained counterparts. These variations present important practical challenges and fairness considerations that must be addressed to design effective and equitable mKEPs.

Table~\ref{tab:lit_comparison} situates this study against the closest related work along the
dimensions discussed above. No prior study combines a direct individual-rationality-versus-axiomatic-fairness comparison with empirical calibration to an Indian
registry, and none simulates blood-group composition asymmetry across registries as a
controlled experimental factor: \citet{Druzsin2024} comes closest, calibrating arrival volumes to
real, country-specific data, but the resulting heterogeneity across countries is a byproduct of
using realistic national instances rather than a factor varied specifically to isolate its effect
on mechanism outcomes, and the underlying comparison is across collaboration policies and
optimization criteria rather than across an IR-anchored and a Shapley-anchored allocation
mechanism.

\begin{table}[ht]
\centering
\footnotesize
\renewcommand{\arraystretch}{1.3}
\resizebox{\textwidth}{!}{%
\begin{tabular}{p{2.6cm}p{3.6cm}ccccc}
\hline
Study & Mechanism & \makecell{Multi-\\period?} & \makecell{Heterog.\\registries?} & \makecell{BG heterog.\\across registries?} & \makecell{Non-Western\\data?} & \makecell{IR vs.\ Shapley\\comparison?} \\
\hline
\citet{Hajaj2015} & Dynamic credit (strategy-proof, long-run IR) & \checkmark & \checkmark (centers) & -- & -- & -- \\
\citet{Biro2020} & Shapley- vs.\ benefit-value-anchored credit & \checkmark & \checkmark & -- & -- & -- \\
\citet{Klimentova2021} & Multi-round credit, several cooperative-game concepts compared & \checkmark & \checkmark & -- & -- & -- \\
\citet{Benedek2024} & Credit-based, deviation-minimizing (large scale, up to 15 countries) & \checkmark & \checkmark & -- & -- & -- \\
\citet{Druzsin2024} & Simulation across collaboration policies and optimization criteria & \checkmark & \checkmark & (\checkmark)\footnotemark & -- & -- \\
\textbf{This paper} & IR-anchored vs.\ Shapley-anchored credit, directly compared & \checkmark & \checkmark & \checkmark & \checkmark (India) & \checkmark \\
\hline
\end{tabular}%
}
\caption{Comparison with the closest related work on multi-agent kidney exchange fairness and
individual-rationality mechanisms. ``Heterog.\ registries'' indicates whether registries/countries
differ in any studied parameter (e.g.\ size, arrival rate); ``BG heterog.\ across registries''
indicates specifically whether blood-group composition asymmetry across registries is examined as
a factor; ``IR vs.\ Shapley comparison'' indicates whether an individual-rationality-anchored
mechanism and a Shapley-value-anchored fairness mechanism are compared directly within the same
study, holding the underlying solution method fixed.}
\label{tab:lit_comparison}
\end{table}
\footnotetext{\citet{Druzsin2024} calibrate arrival volumes and other parameters to real,
country-specific data, which plausibly differs in blood-group composition across the countries
studied; however, blood-group asymmetry is not isolated as a controlled factor manipulated to
study its own effect, as it is in the present study.}

This study presents a comparative simulation analysis of four kidney exchange allocation strategies over a multi-year period: (i) standalone allocation within individual registries; (ii) multi-registry allocation without fairness or individual-rationality constraints, intended to capture the maximum potential benefit of joint operations; (iii) multi-registry allocation with a cumulative individual-rationality mechanism that never revises a registry's surplus downward; and (iv) multi-registry allocation incorporating Shapley-value-based fairness criteria under a signed, budget-balancing rolling credit that redistributes surplus in later rounds. By holding the underlying lexicographic solution method fixed across (iii) and (iv) and varying only how the round's target is constructed and whether surplus is retained or redistributed, we isolate the effect of the fairness mechanism's design from the effect of registry heterogeneity itself, directly addressing the first gap identified above. The analysis focuses on a two-registry setup based on data distributions from Indian kidney exchange registries, an empirical context essentially absent from the existing mKEP fairness literature and directly addressing the second gap. A six-year simulation study is conducted, with matching runs occurring every three months, resulting in a total of 24 simulation rounds; the main findings are further validated at a ten-year (40-round) horizon as a robustness check, confirming that the cumulative individual-rationality and fairness dynamics reported below have stabilized rather than reflecting an artifact of the chosen horizon length. The process is replicated multiple times to account for stochastic variation and to ensure the robustness of the results.

Beyond comparing mechanisms in isolation, we examine how the magnitude of registry heterogeneity in arrival rate, blood-group composition, cycle-length bound, and dropout probability, affects the size and distribution of the achievable gain from collaboration under each mechanism, and discuss the conditions under which that gain may be too small, or too asymmetrically distributed toward one registry, for participation in an mKEP to be worthwhile even under a formal fairness or individual-rationality. Match quality is assessed via two indicators, HLA compatibility and donor--recipient age difference, with a further stratified analysis by blood group to examine differential impacts across blood types.

Across this analysis, three findings stand out. First, pooling increases transplant volume in almost every scenario examined, but this is not unconditional: under blood-group-composition asymmetry specifically, the registry with the more favorable donor composition can realize \emph{fewer} transplants under unconstrained pooling than it would achieve alone, a loss that both the individual-rationality and Shapley-fair mechanisms eliminate. Second, the benefit of pooling is systematically uneven, the more constrained registry, whether by arrival rate, blood-group composition, cycle-length bound, or dropout probability, consistently gains more than its partner, to the point that a registry already well served on its own may find the gain from joining an mKEP too small to justify the cost of participation. Third, a blood-group-stratified breakdown traces the composition-driven loss to a single blood type: O-type recipients bear almost the entire burden of the harm caused by unconstrained pooling, and remain somewhat disadvantaged even under the fairness-oriented mechanisms in several scenarios, though far less severely than without them.

The remainder of the paper is organized as follows. Section \ref{sec: Model} presents integer programming models for the four allocation strategies. Section \ref{sec: Data} outlines the data and simulation methodology employed in the study. The results of the analysis are presented in Section \ref{sec: Results}, followed by a concluding discussion in Section \ref{sec: Conclusion}.

\section{Models and allocation process} \label{sec: Model}

We compare four allocation mechanisms for a multi-registry kidney exchange program (mKEP):
(i) independent, registry-level allocation; (ii) pooled mKEP allocation without any fairness or
individual-rationality (IR) constraints; (iii) pooled mKEP allocation with an individual-rationality
constraint, tracked cumulatively across matching rounds; and (iv) pooled mKEP allocation with
fairness enforced through Shapley-value-based compensation. We formalize each mechanism as an
integer program (IP), building on the compact edge formulation of \citet{Constantino2013}.

\subsection{Notation and common formulation} \label{sec: notation}

We represent the kidney exchange problem as a directed graph $G=(V,A)$, where each node in $V$
corresponds to a donor--recipient pair and each arc $(i,j)\in A$ indicates that pair $i$'s donor
is compatible with pair $j$'s recipient. Let $R=\{1,\dots,|R|\}$ index the $|R|$ participating
registries, each with its own disjoint node set $V_k$, so that $V=\bigcup_{k\in R}V_k$.

We replicate $G$ a total of $M=\lceil |V|/2\rceil$ times, an
upper bound on the number of vertex-disjoint cycles the compatibility graph can support, in order
to bound individual exchange-cycle size without resorting to an exponential cycle-enumeration
formulation. For each replication $l\in\{1,\dots,M\}$ and arc $(i,j)\in A$, define the binary
decision variable
\begin{equation}
  x_{ij}^{l} =
  \begin{cases}
    1 & \text{if arc } (i,j) \text{ is selected in the } l\text{-th replication of } G, \\
    0 & \text{otherwise},
  \end{cases}
  \label{mKEP variable}
\end{equation}
which, with the node set restricted appropriately, is used in all four models below. Let
$w_{ij}^{l}\ge0$ denote the match-quality weight of arc $(i,j)$ in replication $l$.

Each registry $k$ imposes its own intra-registry cycle-length bound $b_k$, collected in the
vector $b=(b_1,\dots,b_{|R|})$; once registries pool their pairs, a common cycle-length bound $B$
applies to every replication. These bounds reflect institutional, logistical, or policy
limitations on the size of a permissible exchange.

For a node subset $V'\subseteq V$, define the feasible region $\mathcal{F}(V')$ by
\begin{align}
  \sum_{j:(j,i)\in A} x_{ji}^{l} &= \sum_{j:(i,j)\in A} x_{ij}^{l},
    && \forall i \in V', \ \forall l \in \{1,\dots,M\}
    \label{eq:flowbalance}\\
  \sum_{l=1}^{M}\sum_{j:(i,j)\in A} x_{ij}^{l} &\le 1,
    && \forall i \in V'
    \label{eq:max_donation}\\
  x_{ij}^{l} &\in \{0,1\},
    && \forall (i,j)\in A,\ \forall l \in \{1,\dots,M\}.
    \label{eq:binary}
\end{align}
Constraint~\eqref{eq:flowbalance} enforces simultaneity of exchange within a replication: a
pair's donor contributes a kidney only if that pair's recipient also receives one, in the same
replication. Constraint~\eqref{eq:max_donation} restricts each pair to at most one transplant,
summed across all replications. Constraint~\eqref{eq:binary} imposes integrality. We reuse
$\mathcal{F}(\cdot)$ throughout: the individual model applies it to a single registry's node set
$V_k$, while all three mKEP models apply it to the pooled set $V$, differing only in objective
and in any fairness-related constraints layered on top.

The registry-level and cross-registry bounds introduced in Sections~\ref{sec: Individual
allocation}--\ref{sec: mKEP_without_fair} are each stated as a single linear inequality summed
over an entire replica. The following lemma establishes that, despite this aggregated form, each
such bound is an \emph{exact} restriction on every individual cycle realized in that replica, not
merely a cap on their combined length.

\begin{lemma} \label{lem:cycle_decomp}
Fix a replica $l\in\{1,\dots,M\}$ and let $x\in\mathcal{F}(V')$ for some $V'\subseteq V$. Then:
\begin{enumerate}
\item[(a)] The arcs $\{(i,j)\in A : x_{ij}^{l}=1\}$ decompose into vertex-disjoint simple
directed cycles $C_1,\dots,C_p$ within replica $l$ (possibly $p=0$).
\item[(b)] For any subset $W\subseteq V'$ and any $b\ge0$, if
$\sum_{(i,j)\in A:\,i\in W} x_{ij}^{l} \le b$, then $|C_m\cap W|\le b$ for every $m=1,\dots,p$.
\end{enumerate}
\end{lemma}
\begin{proof}
(a) Constraint~\eqref{eq:max_donation}, summed over all $M$ replicas, restricts each $i\in V'$ to
at most one selected outgoing arc in total; since $x$ is binary, node $i$ therefore has out-degree
$1$ in at most one replica and out-degree $0$ in every other. Fix replica $l$ and call a node
\emph{active} in $l$ if its out-degree there is $1$. By constraint~\eqref{eq:flowbalance},
in-degree equals out-degree for every $i\in V'$ within replica $l$; hence every active node also
has in-degree exactly $1$ in $l$, and every inactive node has in-degree $0$. Restricted to the
active nodes, replica $l$'s selected arcs therefore form a functional digraph with in-degree $=$
out-degree $=1$ at every node. Following any active node forward along its unique outgoing arc
must, after finitely many steps, revisit a node; the in-degree-$1$ property forces that revisited
node to be the starting node, closing a simple cycle. Two such cycles cannot share a node, since a
shared node would require out-degree at least $2$. This gives the claimed decomposition into
vertex-disjoint simple cycles $C_1,\dots,C_p$.

(b) Since $C_1,\dots,C_p$ are vertex-disjoint and each node of $C_m\cap W$ contributes exactly one
arc to the sum (its unique outgoing arc within $C_m$), the sum decomposes as
\[
\sum_{(i,j)\in A:\,i\in W} x_{ij}^{l} \;=\; \sum_{m=1}^{p} |C_m\cap W|.
\]
By hypothesis this total is at most $b$. Each term $|C_m\cap W|$ is a non-negative integer, and a
non-negative summand cannot exceed a sum of non-negative summands bounded by $b$; hence
$|C_m\cap W|\le b$ for every $m$.
\end{proof}

\begin{remark}
Applying Lemma~\ref{lem:cycle_decomp}(b) with $W=V_k$ and $b=b_k$ shows that the registry-level
bound in~\eqref{eq:registry_bound_single} and~\eqref{eq:registry_bound_mKEP} restricts registry
$k$'s participation in \emph{every} cycle realized in a given replica, not merely their combined
length across the replica's (possibly several) cycles. Applying it with $W=V$ and $b=B$ gives the
analogous exact restriction for the cross-registry bound in~\eqref{eq:global_bound_mKEP} on total
cycle length. This is the basis for the claim in Section~\ref{sec: Data} that a bound applied
within a single replica is an exact -- not approximate -- restriction on the length of every
cycle realized in that replica, regardless of how many separate cycles the replica's arc budget
happens to be split across.
\end{remark}

\subsection{Individual allocation model} \label{sec: Individual allocation}

The individual allocation model represents registry $k$ operating alone. It maximizes the
weighted number of transplants realized among $k$'s own pairs, subject to $k$'s own
cycle-length bound:
\begin{align}
  \max \quad & \sum_{l=1}^{M} \sum_{(i,j)\in A:\, i\in V_k} w_{ij}^{l}\, x_{ij}^{l}
    \label{eq:objective_single}\\
  \text{s.t.} \quad & x \in \mathcal{F}(V_k) \label{eq:feas_single}\\
  & \sum_{(i,j)\in A:\, i\in V_k} x_{ij}^{l} \le b_k, \qquad \forall l\in\{1,\dots,M\}.
    \label{eq:registry_bound_single}
\end{align}
We denote the \emph{number} of transplants (not the weighted score) obtained by this model for
registry $k$ at round $t$ by $IR_k^{(t)}$; this quantity is registry $k$'s standalone benchmark
in Section~\ref{sec: mkep_with_IR}.

\subsection{mKEP model without fairness or individual-rationality constraints}
\label{sec: mKEP_without_fair}

Registries may agree to pool their pairs into a single compatibility graph without any explicit
safeguard for how the resulting gains are distributed. This model maximizes system-wide weighted
transplant volume:
\begin{align}
  \max \quad & \sum_{l=1}^{M}\sum_{(i,j)\in A} w_{ij}^{l}\, x_{ij}^{l}
    \label{eq:objective_mKEP} \\[1ex]
  \text{s.t.} \quad
  & x \in \mathcal{F}(V) \label{eq:feas_mKEP} \\[1ex]
  & \sum_{j:(i,j)\in A} x_{ij}^{l} \leq b_{k},
    && \forall k \in R, \, \forall i \in V_k, \, \forall l \in \{1,\dots,M\}
    \label{eq:registry_bound_mKEP} \\[1ex]
  & \sum_{(i,j) \in A} x_{ij}^{l} \leq B,
    && \forall l \in \{1,\dots,M\}.
    \label{eq:global_bound_mKEP}
\end{align}
This model represents full data-sharing with no correction for how gains are distributed: a
registry may realize fewer transplants within the mKEP than it would have achieved on its own.
Sections~\ref{sec: mkep_with_IR} and~\ref{sec: mkep_with_fair} address this through two
complementary mechanisms -- a minimal, count-based individual-rationality guarantee, and a
Shapley-value-based fair-share target -- both formulated so as to remain feasible in every round,
regardless of how large a shortfall may have accumulated in earlier rounds.

\subsection{mKEP model with individual rationality constraints} \label{sec: mkep_with_IR}

A minimal requirement for a registry's continued participation in the mKEP is that it should
never, in expectation, be made worse off than by acting alone. We formalize this as an
\emph{individual rationality} (IR) requirement, applied both within a single round and
cumulatively over the matching horizon.

\subsubsection{Round-level and cumulative individual rationality}

Let $A_k^{(t)}=\sum_{l=1}^{M}\sum_{(i,j)\in A:\, i\in V_k} x_{ij}^{l,(t)}$ denote registry $k$'s
realized number of transplants within the mKEP at round $t$, and let $IR_k^{(t)}$ be its
standalone benchmark from Section~\ref{sec: Individual allocation}, computed on the pairs present
at round $t$. Round-level IR requires $A_k^{(t)}\ge IR_k^{(t)}$ for every $k\in R$.

\begin{proposition} \label{prop:feasibility}
For every round $t$, there exists a feasible mKEP solution satisfying $A_k^{(t)}\ge IR_k^{(t)}$
simultaneously for all $k\in R$.
\end{proposition}
\begin{proof}
Let $x^{k,\star}$ be an optimal solution to the individual allocation model
(Section~\ref{sec: Individual allocation}) for registry $k$'s pairs present at round $t$, using
only arcs internal to $V_k$. Because the sets $V_k$ are pairwise disjoint and each $x^{k,\star}$
satisfies~\eqref{eq:flowbalance}--\eqref{eq:binary} using only nodes in $V_k$, the union
$\bigcup_{k\in R} x^{k,\star}$ satisfies~\eqref{eq:flowbalance}--\eqref{eq:binary} on $V$ and
achieves exactly $IR_k^{(t)}$ transplants for every registry $k$ simultaneously. This solution is
feasible for the mKEP model of Section~\ref{sec: mKEP_without_fair}, which can only do at least
as well by additionally exploiting cross-registry arcs.
\end{proof}

Registries are also concerned with their \emph{cumulative} standing, not merely round-level
outcomes. Let $C_k^{IS}(t)=\sum_{t'=1}^{t} IR_k^{(t')}$ and
$C_k^{mKEP}(t)=\sum_{t'=1}^{t} A_k^{(t')}$ denote registry $k$'s cumulative independent and
cumulative mKEP totals through round $t$, and define the carried-forward deficit
\begin{equation}
  D_k(t) = \max\bigl(0,\ C_k^{IS}(t) - C_k^{mKEP}(t)\bigr), \qquad D_k(0)=0.
  \label{eq:deficit}
\end{equation}
Unlike Proposition~\ref{prop:feasibility}'s round-level guarantee, a \emph{hard} constraint
enforcing $D_k(t)=0$ for all $t$ can become infeasible: the mKEP and the independent benchmark are
distinct dynamical systems that, having matched different pairs in earlier rounds, retain
different residual pools, so nothing guarantees the mKEP can recoup an accumulated deficit within
the capacity of the remaining rounds. We therefore treat cumulative IR as a soft, lexicographically
prioritized objective rather than a hard constraint.

\subsubsection{Lexicographic formulation}

At round $t$, define the target
\begin{equation}
  \tau_k^{(t)} = IR_k^{(t)} + D_k(t-1), \qquad \forall k \in R,
  \label{eq:IRtarget}
\end{equation}
and a nonnegative, one-sided deviation variable $y_k$ capturing only shortfall below this target
(a registry realizing more than $\tau_k^{(t)}$ is never penalized, since surplus can only reduce
its future deficit):
\begin{equation}
  y_k \ge \tau_k^{(t)} - A_k^{(t)}, \qquad y_k \ge 0, \qquad \forall k \in R.
  \label{eq:IRdeviation}
\end{equation}
The model is solved in two lexicographically ordered stages. Stage 1 minimizes the aggregate
individual-rationality shortfall:
\begin{align}
  \gamma^\star = \min \quad & \sum_{k\in R} y_k \label{eq:IR_stage1_obj}\\
  \text{s.t.} \quad & x \in \mathcal{F}(V), \ \eqref{eq:registry_bound_mKEP}, \
    \eqref{eq:global_bound_mKEP}, \ \eqref{eq:IRdeviation}. \label{eq:IR_stage1_feas}
\end{align}
By Proposition~\ref{prop:feasibility}, $\gamma^\star$ is always finite, so Stage 1 -- and hence
Stage 2 -- is always feasible. Stage 2 maximizes weighted transplant volume among all solutions
attaining $\gamma^\star$:
\begin{align}
  \max \quad & \sum_{l=1}^{M}\sum_{(i,j)\in A} w_{ij}^{l}\, x_{ij}^{l} \label{eq:IR_stage2_obj}\\
  \text{s.t.} \quad & x \in \mathcal{F}(V), \ \eqref{eq:registry_bound_mKEP}, \
    \eqref{eq:global_bound_mKEP}, \ \eqref{eq:IRdeviation}, \
    \sum_{k\in R} y_k \le \gamma^\star. \label{eq:IR_stage2_feas}
\end{align}
The deficit is then updated using the round's realized transplants $A_k^{(t),\star}$ at the
Stage-2 optimum,
\begin{equation}
  D_k(t) = \max\bigl(0,\ \tau_k^{(t)} - A_k^{(t),\star}\bigr).
  \label{eq:deficit_update}
\end{equation}
Individual rationality is thus defined solely in terms of transplant \emph{counts}: both
$\tau_k^{(t)}$ and $y_k$ are stated in terms of $A_k^{(t)}$, never the weighted
objective~\eqref{eq:IR_stage2_obj}. Registries value transplant volume -- from which their
revenue and long-run viability derive -- more directly than the relative quality of any single
match, so the IR criterion and the efficiency objective are deliberately kept separate; the same
principle underlies the fairness criterion of Section~\ref{sec: mkep_with_fair}.

\subsection{mKEP model with fair allocation criteria} \label{sec: mkep_with_fair}

One of the central challenges in mKEP design is data sharing, as registries may have a strategic
incentive to conceal certain donor--patient pairs. In particular, a registry might internally
match its most compatible pairs and disclose only its harder-to-match pairs to the shared pool.
Figure~\ref{fig: Incentive_to_drop} illustrates this in a simplified two-registry setting, $R_1$
and $R_2$. Registry $R_1$ holds three pairs ($P_1$, $P_2$, $P_3$); registry $R_2$ holds three
pairs ($P_4$, $P_5$, $P_6$). Compatibility arcs connect each pair in $R_1$ with a corresponding
pair in $R_2$, and arc weights represent match quality. Under the mKEP-optimal solution, each
registry achieves two transplants, with $P_1$ matched to $P_4$ and $P_3$ to $P_6$, yielding the
highest overall match quality. Registry $R_1$ may instead observe that $P_1$ and $P_2$ can be
matched internally; by carrying out this internal exchange and disclosing only $P_3$ to the
shared pool, $R_1$ still secures two transplants for itself, but of lower average quality, and at
the cost of overall system efficiency and of $R_2$'s outcomes. Registry $R_2$ faces a symmetric
incentive to match $P_5$ with $P_6$ internally, disclosing only $P_4$. Such self-interested
behavior is individually rational in isolation but reduces both the aggregate volume and quality
of transplants achievable through collaboration.

\begin{figure}[H]
\centering

\begin{minipage}{0.45\textwidth}
\centering
\begin{tikzpicture}[
    pair/.style={circle, draw, minimum size=1.2cm, thick},
    every edge/.style={draw, thick, -{Latex}, font=\small},
    box/.style={draw, rounded corners, thick, inner sep=0.5cm}
]

\node[pair] (P1) at (0,0) {P1};
\node[pair] (P2) at (2.5,0) {P2};
\node[pair] (P3) at (5,0) {P3};
\node[pair] (P4) at (0,-3) {P4};
\node[pair] (P5) at (2.5,-3) {P5};
\node[pair] (P6) at (5,-3) {P6};

\node[box, fit=(P1) (P2) (P3), label=above:\textbf{Registry R1}] {};
\node[box, fit=(P4) (P5) (P6), label=below:\textbf{Registry R2}] {};

\path (P1) edge[bend left=20] node[above] {1} (P2);
\path (P2) edge[bend left=20] node[below] {1} (P1);
\path (P5) edge[bend left=20] node[above] {1} (P6);
\path (P6) edge[bend left=20] node[below] {1} (P5);

\path (P1) edge[bend left=15] node[right] {2} (P4);
\path (P4) edge[bend left=15] node[left] {2} (P1);
\path (P3) edge[bend left=15] node[right] {2} (P6);
\path (P6) edge[bend left=15] node[left] {2} (P3);

\draw[->, thick] (6.2, -1.2) -- ++(2, 5.8) node[midway, above, sloped] {\textbf{Suboptimal}};
\draw[->, thick] (6.2, -1.2) -- ++(2, 0) node[midway, above, sloped] {\textbf{Optimal}};
\draw[->, thick] (6.2, -1.2) -- ++(2, -5.8) node[midway, below, sloped] {\textbf{Suboptimal}};

\end{tikzpicture}
\end{minipage}
\hfill
\begin{minipage}{0.42\textwidth}
\centering

\begin{tikzpicture}[
    pair/.style={circle, draw, minimum size=1cm, thick},
    match/.style={draw=red!70!black, thick, -{Latex}, line width=1.2pt},
    box/.style={draw, rounded corners, thick, inner sep=0.5cm}
]
\node[pair] (P1) at (0,0) {P1};
\node[pair] (P2) at (2.5,0) {P2};
\node[pair] (P3) at (5,0) {P3};
\node[pair] (P4) at (0,-3) {P4};
\node[pair] (P5) at (2.5,-3) {P5};
\node[pair] (P6) at (5,-3) {P6};
\node[box, fit=(P1) (P2) (P3), label=above:\textbf{Registry 1 Biased Solution}] {};
\node[box, fit=(P4) (P5) (P6)] {};
\path[match] (P1) edge[<->, >=Latex, left=15] node[above] {2} (P2);
\path[match] (P3) edge[<->, >=Latex, left=15] node[right] {4} (P6);
\end{tikzpicture}

\vspace{0.5cm}

\begin{tikzpicture}[
    pair/.style={circle, draw, minimum size=1cm, thick},
    match/.style={draw=green!60!black, thick, -{Latex}, line width=1.2pt},
    box/.style={draw, rounded corners, thick, inner sep=0.5cm}
]
\node[pair] (P1) at (0,0) {P1};
\node[pair] (P2) at (2.5,0) {P2};
\node[pair] (P3) at (5,0) {P3};
\node[pair] (P4) at (0,-3) {P4};
\node[pair] (P5) at (2.5,-3) {P5};
\node[pair] (P6) at (5,-3) {P6};
\node[box, fit=(P1) (P2) (P3), label=above:\textbf{Optimal Solution}] {};
\node[box, fit=(P4) (P5) (P6)] {};
\path[match] (P1) edge[<->, >=Latex, left=15] node[right] {4} (P4);
\path[match] (P3) edge[<->, >=Latex, left=15] node[right] {4} (P6);
\end{tikzpicture}

\vspace{0.5cm}

\begin{tikzpicture}[
    pair/.style={circle, draw, minimum size=1cm, thick},
    match/.style={draw=yellow!70!black, thick, -{Latex}, line width=1.2pt},
    box/.style={draw, rounded corners, thick, inner sep=0.5cm}
]
\node[pair] (P1) at (0,0) {P1};
\node[pair] (P2) at (2.5,0) {P2};
\node[pair] (P3) at (5,0) {P3};
\node[pair] (P4) at (0,-3) {P4};
\node[pair] (P5) at (2.5,-3) {P5};
\node[pair] (P6) at (5,-3) {P6};
\node[box, fit=(P1) (P2) (P3), label=above:\textbf{Registry 2 Biased Solution}] {};
\node[box, fit=(P4) (P5) (P6)] {};
\path[match] (P1) edge[<->, >=Latex, left=15] node[right] {4} (P4);
\path[match] (P5) edge[<->, >=Latex, left=15] node[above] {2} (P6);
\end{tikzpicture}

\end{minipage}
\caption{Original kidney exchange graph (left) with arrows pointing to optimal and suboptimal
matchings (right)}
\label{fig: Incentive_to_drop}
\end{figure}

To counteract this incentive, we allocate transplants toward a registry-specific \emph{target
solution}, $TS_k$, based on the Shapley value of the pooled coalition together with a rolling
credit for past shortfalls, and solve for the mKEP allocation that comes closest to these
targets. As in Section~\ref{sec: mkep_with_IR}, fairness is defined solely in terms of transplant
\emph{counts}: incorporating match-quality weights into the fairness criterion itself would risk
undervaluing the credit owed to a registry for transplants missed in earlier rounds, since one
additional transplant is of greater value to a registry's throughput and revenue than the
marginal quality difference between any two matches.

\subsubsection{Shapley fair-share target}

We use the Shapley value specifically, rather than an alternative cooperative-game solution
concept, following evidence that it outperforms the alternatives on this task: \citet{Klimentova2021}
and \citet{Biro2020} find it produces more equitable outcomes than the benefit or contribution
value in multi-round kidney exchange simulations, and \citet{Benedek2024} independently confirm
this at considerably larger scale, finding a Shapley-value-based target outperforms the nucleolus,
Banzhaf value, and tau value alike in large international kidney exchange schemes.

For a coalition $S\subseteq R$, let $V_S=\bigcup_{k\in S}V_k$ and define the coalition's
standalone value as the optimal \emph{count}-based objective of the mKEP model of
Section~\ref{sec: mKEP_without_fair}, restricted to $V_S$:
\begin{equation}
  v(S) = \max \Bigl\{ \textstyle\sum_{l=1}^{M}\sum_{(i,j)\in A:\, i,j\in V_S} x_{ij}^{l}
  \ :\ x\in\mathcal{F}(V_S),\ \eqref{eq:registry_bound_mKEP}\ \forall k\in S,\
  \eqref{eq:global_bound_mKEP} \Bigr\}.
  \label{eq:coalition_value}
\end{equation}
Registry $k$'s Shapley fair share of the grand coalition $R$ is
\begin{equation}
  \phi_k = \sum_{S\subseteq R\setminus\{k\}} \frac{|S|!\,(|R|-|S|-1)!}{|R|!}\,
  \bigl[v(S\cup\{k\}) - v(S)\bigr],
  \label{eq:shapley}
\end{equation}
which, for the two-registry case $R=\{1,2\}$ used in our computational study, reduces to
\begin{equation}
  \phi_1 = v(\{1\}) + \tfrac12\bigl(v(R) - v(\{1\}) - v(\{2\})\bigr), \qquad
  \phi_2 = v(\{2\}) + \tfrac12\bigl(v(R) - v(\{1\}) - v(\{2\})\bigr).
  \label{eq:shapley2}
\end{equation}

\subsubsection{Rolling credit and lexicographic formulation}

Analogously to the individual-rationality deficit of Section~\ref{sec: mkep_with_IR}, let
$\rho_k(t-1)\in\mathbb{R}$ denote registry $k$'s carried-forward fairness credit entering round
$t$ ($\rho_k(0)=0$; unlike the IR deficit $D_k(t-1)$, $\rho_k(t-1)$ is \emph{signed} rather than
truncated at zero -- see below), and define the round's target solution
\begin{equation}
  TS_k^{(t)} = \phi_k^{(t)} + \rho_k(t-1).
  \label{eq:TS}
\end{equation}
Whereas the individual-rationality mechanism of Section~\ref{sec: mkep_with_IR} penalizes only
shortfall below target, the fairness mechanism here penalizes deviation in \emph{either}
direction: a registry realizing more than its target $TS_k^{(t)}$ is also penalized, since its
surplus is deducted from next round's target rather than retained unconditionally. This is a
deliberate point of contrast between the two mechanisms compared in this study: the IR mechanism
never claws back a registry's surplus, whereas the Shapley-based mechanism here always does,
which (subject to Remark~\ref{rem:budget}) keeps it budget-balanced.
\begin{equation}
  y_k \ge TS_k^{(t)} - A_k^{(t)}, \qquad y_k \ge A_k^{(t)} - TS_k^{(t)}, \qquad y_k \ge 0,
  \qquad \forall k\in R.
  \label{eq:fair_deviation}
\end{equation}
Equivalently, $y_k \ge \bigl|TS_k^{(t)} - A_k^{(t)}\bigr|$. (A linear objective of the form
$\min\sum_k \bigl(A_k^{(t)} - TS_k^{(t)}\bigr)$, without deviation variables, is unbounded below
in $x$ and does not enforce proximity to the target regardless of which penalty structure is
chosen; \eqref{eq:fair_deviation} is what corrects this.) The model is solved in the same two
lexicographically ordered stages as Section~\ref{sec: mkep_with_IR}:
\begin{align}
  \gamma^\star = \min \quad & \sum_{k\in R} y_k
    && \text{s.t.} \ x\in\mathcal{F}(V),\ \eqref{eq:registry_bound_mKEP},\
    \eqref{eq:global_bound_mKEP},\ \eqref{eq:fair_deviation}
    \label{eq:fair_stage1}\\
  \max \quad & \sum_{l=1}^{M}\sum_{(i,j)\in A} w_{ij}^{l}\, x_{ij}^{l}
    && \text{s.t.} \ x\in\mathcal{F}(V),\ \eqref{eq:registry_bound_mKEP},\
    \eqref{eq:global_bound_mKEP},\ \eqref{eq:fair_deviation},\
    \textstyle\sum_{k\in R} y_k \le \gamma^\star
    \label{eq:fair_stage2}
\end{align}
with the credit updated, \emph{without} truncation at zero, as
\begin{equation}
  \rho_k(t) = TS_k^{(t)} - A_k^{(t),\star},
  \label{eq:credit_update}
\end{equation}
using the transplant \emph{count} realized at the Stage-2 optimum. Feasibility of Stage~1 follows
by the same argument as Proposition~\ref{prop:feasibility}: since $v(\{k\})$ is achievable by
registry $k$ acting alone, the union of each registry's standalone-optimal solution remains a
feasible point for the pooled model regardless of the sign of the resulting deviation, which
bounds $\gamma^\star$ and guarantees Stage~1 -- and hence Stage~2 -- is always solvable, however
large $|\rho_k(t-1)|$ has grown.

\begin{remark}[Budget balance] \label{rem:budget}
Because the Shapley value satisfies efficiency, $\sum_{k\in R}\phi_k^{(t)} = v(R)^{(t)}$ in every
round, where $v(R)^{(t)}$ is the coalition value~\eqref{eq:coalition_value} evaluated on the pairs
present at round $t$. If, in addition, the Stage-2-optimal solution attains total transplant
count exactly $v(R)^{(t)}$ -- i.e., no volume is sacrificed in resolving the Stage-1 fairness
objective -- then $\sum_{k\in R}\rho_k(t) = 0$ and, inductively, $\sum_{k\in R}TS_k^{(t+1)} =
v(R)^{(t+1)}$: the credit mechanism exactly redistributes the coalition's realized value among
registries, without gain or loss. This coincidence holds in the worked example below, but is not
guaranteed in general under the fairness-first lexicographic ordering
of~\eqref{eq:fair_stage1}--\eqref{eq:fair_stage2}: Stage~1 may, in principle, select a
lower-volume solution if doing so strictly reduces $\sum_k y_k$, in which case
$\sum_k \rho_k(t) < 0$ at that round specifically, even though the mechanism remains feasible and
well defined.
\end{remark}

\subsubsection{Illustrative example}

Table~\ref{tab:credit_example} traces the mechanism over two rounds for a three-registry
illustration -- a generalization beyond the two-registry case used elsewhere in this study,
included to show that the mechanism applies unchanged for $|R|>2$ via the general Shapley
formula~\eqref{eq:shapley}. Initial credit is zero for every registry. In round $t=1$, the
computed Shapley shares are $\phi^{(1)}=(6.8,\,22.6,\,10.6)$; with no prior credit, the target
solution equals the Shapley share exactly. The closest attainable matching realizes
$A^{(1),\star}=(8.0,\,22.0,\,10.0)$: registry~1 exceeds its target while registries~2 and~3 fall
short, yielding credits $\rho^{(1)}=(-1.2,\,0.6,\,0.6)$ via~\eqref{eq:credit_update}. In round
$t=2$, freshly computed Shapley shares $\phi^{(2)}=(10.2,\,12.8,\,7.0)$ combine with the carried
credit to give the target $TS^{(2)}=(9.0,\,13.4,\,7.6)$: registry~1's earlier surplus is
subtracted from its new target, while registries~2 and~3's earlier shortfalls are added to
theirs. The closest attainable matching this round is $A^{(2),\star}=(8.0,\,14.0,\,8.0)$, giving
updated credits $\rho^{(2)}=(1.0,\,-0.6,\,-0.4)$.

\begin{table}[H]
\centering
\renewcommand{\arraystretch}{1.3}
\begin{minipage}[t]{0.48\textwidth}
\centering
\textbf{Round $t=1$}\\[2pt]
\resizebox{\linewidth}{!}{%
\begin{tabular}{lcccc}
\toprule
 & $R_1$ & $R_2$ & $R_3$ & Total \\
\midrule
$\phi^{(1)}$ & 6.8 & 22.6 & 10.6 & 40.0 \\
$TS^{(1)}=\phi^{(1)}+\rho(0)$ & 6.8 & 22.6 & 10.6 & 40.0 \\
$A^{(1),\star}$ & 8.0 & 22.0 & 10.0 & 40.0 \\
$\rho^{(1)}=TS^{(1)}-A^{(1),\star}$ & $-1.2$ & 0.6 & 0.6 & 0.0 \\
\bottomrule
\end{tabular}%
}
\end{minipage}%
\hfill
\begin{minipage}[t]{0.48\textwidth}
\centering
\textbf{Round $t=2$}\\[2pt]
\resizebox{\linewidth}{!}{%
\begin{tabular}{lcccc}
\toprule
 & $R_1$ & $R_2$ & $R_3$ & Total \\
\midrule
$\phi^{(2)}$ & 10.2 & 12.8 & 7.0 & 30.0 \\
$TS^{(2)}=\phi^{(2)}+\rho^{(1)}$ & 9.0 & 13.4 & 7.6 & 30.0 \\
$A^{(2),\star}$ & 8.0 & 14.0 & 8.0 & 30.0 \\
$\rho^{(2)}=TS^{(2)}-A^{(2),\star}$ & 1.0 & $-0.6$ & $-0.4$ & 0.0 \\
\bottomrule
\end{tabular}%
}
\end{minipage}
\caption{Two-sided credit mechanism traced over two rounds for a three-registry illustration.
Registry~1's round-1 surplus (achieving 8.0 against a target of 6.8) is deducted from its
round-2 target; registries~2 and~3's round-1 shortfalls are added to theirs.}
\label{tab:credit_example}
\end{table}
 
In both rounds the realized total equals the target sum exactly (40.0 and 30.0, respectively), so
this illustration falls in the budget-balanced case of Remark~\ref{rem:budget}: here the credit
mechanism purely reallocates a fixed pool of realized value among registries, compensating a
registry that received less than its fair share in the current round by prioritizing it in the
next.

\subsection{Summary of the four models}

Table~\ref{tab:summary} summarizes the four models. Models (iii) and (iv) share the same overall
architecture - a lexicographic solve prioritizing closeness to a rolling target ahead of weighted
volume but differ along two independent dimensions, each a deliberate point of comparison in
this study. First, the target itself is constructed differently: model (iii) anchors to registry
$k$'s own realized standalone trajectory, while model (iv) anchors to its axiomatic (Shapley)
share of the coalition's value. Second, the two mechanisms penalize deviation from target
differently: model (iii) penalizes shortfall only, so a registry's surplus is retained
unconditionally and never reduces a future target, whereas model (iv) penalizes deviation in
either direction, clawing back surplus via the signed rolling credit
of~\eqref{eq:credit_update} and keeping the mechanism budget-balanced under the conditions of
Remark~\ref{rem:budget}. Comparing (iii) against (iv) therefore isolates the effect of the credit
rule's symmetry in particular, whether registries with easier-to-match patient pools retain
their surplus indefinitely (iii) or have it redistributed to their partner registries (iv) from
the separate question of how the target itself is constructed. Both mechanisms remain feasible in
every round by the same underlying argument (Proposition~\ref{prop:feasibility}), which is what
permits treating a registry's fairness or IR shortfall as a reported outcome of the model rather
than a source of infeasibility.

\begin{table}[H]
\centering
\begin{tabular}{p{3.2cm}p{4.3cm}p{6.3cm}}
\hline
Model & Feasible region & Objective \\
\hline
(i) Individual & $\mathcal{F}(V_k)$, per registry & Maximize weighted transplants, registry $k$ \\
(ii) mKEP, no fairness & $\mathcal{F}(V)$ & Maximize weighted transplants, system-wide \\
(iii) mKEP, IR & $\mathcal{F}(V)$, target $\tau_k^{(t)}=IR_k^{(t)}+D_k(t-1)$,\ $D_k\ge0$ &
  Lexicographic: minimize \emph{one-sided} count shortfall, then maximize weighted transplants \\
(iv) mKEP, Shapley-fair & $\mathcal{F}(V)$, target $TS_k^{(t)}=\phi_k^{(t)}+\rho_k(t-1)$,\ $\rho_k\in\mathbb{R}$ &
  Lexicographic: minimize \emph{two-sided} count deviation, then maximize weighted transplants \\
\hline
\end{tabular}
\caption{Summary of the four allocation models.}
\label{tab:summary}
\end{table}

\section{Data Sources and Simulation Design} \label{sec: Data}

The study uses data derived from the Apex Swap Transplant Registry (ASTRA), an Indian registry
dedicated to swap transplants \citep{Billa2018}, together with a second blood-group distribution
derived from Roth et al.'s reported Alliance for Paired Kidney Donation (APKD) data
\citep{Ashlagi2021}. Assigning ASTRA's distribution to one registry and the APKD-derived
distribution to the other allows the two simulated registries to differ systematically in
composition rather than share an identical blood-group mix. The aim is to evaluate the benefits
conferred by a multi-registry framework in terms of both transplant volume and match quality.
Transplant compatibility is determined by blood-group matching; match quality is assessed via two
indicators: the number of HLA matches (the scoring scheme is detailed in
\ref{sec: Scoring_system}) and the age difference between donor and recipient \citep{Billa2018}.

We examines the mKEP framework along four dimensions: (a) the incremental benefit accruing
to each registry by participating in the mKEP
under each of the four allocation mechanisms described in Section~\ref{sec: Model}; (b) the
influence of registry size, operationalized as arrival rate (c) the effect of differences in blood-group composition across registries;
and (d) the effect of differences in the per-registry cycle-length bound and in dropout
probability across registries.

A two-registry configuration is used for the simulation, though the underlying models extend
directly to more than two registries, as illustrated by the general Shapley
formula~\eqref{eq:shapley} and the three-registry example of Section~\ref{sec: mkep_with_fair};
\citet{Benedek2024} demonstrate computationally that a closely related credit-based,
deviation-minimizing mechanism remains tractable at considerably larger scale, up to fifteen
countries. The simulation spans 24 matching rounds, corresponding to
four rounds per year over six years, and is replicated 50 times per configuration to average out
the stochastic variation in blood-group generation, arrival counts, and dropout outcomes; a
10-year (40-round) horizon is additionally simulated as a robustness check \ref{sec: Appendix_Robustness}, to confirm that the
cumulative individual-rationality and fairness dynamics reported at 24 rounds have stabilized
rather than reflecting an artifact of the chosen horizon length.

All integer programs were solved using Gurobi Optimizer version 12.0.3,
accessed through its Python interface (\texttt{gurobipy}), under an academic license. Simulations
were run on an Apple MacBook Pro (M3 chip, macOS) with 8~GB of RAM and 512~GB of SSD storage.
Solve time for a complete simulation run with 50 replications of 24 rounds each, across all four
allocation models ranged from a few minutes to approximately an hour, depending on the
scenario configuration.

Arrivals are drawn from a Poisson process with mean 5 (Small) or mean 10 (Medium) pairs per round
for the two registries; as a robustness check on this modeling choice, arrivals are additionally
drawn from a Negative Binomial process with matched mean and an explicit over dispersion parameter,
to confirm that the reported findings are not an artifact of assuming Poisson (equidispersed)
arrivals specifically. Dropout probability -- the likelihood that an unmatched pair exits the
registry before the following round is set to Low (0.2) or High (0.4), and a pair can stay in
the registry for a maximum of 2 years (i.e., 8 rounds) \citep{Verissimo2022}.

Each registry $k$ is subject to a cycle-length bound $b_k$: within any one matching round, no
single exchange cycle in which registry $k$ participates may involve more than $b_k$ of registry
$k$'s pairs. This is not a cap on registry $k$'s total round volume -- a registry may still
complete several separate cycles within the same three-month matching period, each individually
respecting the bound -- but a structural limit on how large any one coordinated exchange is
permitted to be, reflecting the logistical burden of arranging simultaneous donor nephrectomies
across a growing number of participating pairs. We compare a restrictive setting ($b_k=2$, restricting
each registry to simple two-way swaps internally) against a relaxed setting ($b_k=3$), and, for
cross-registry cycles specifically, a separate bound $B$ ranging from 3 to 5, reflecting the
additional coordination overhead of a cycle spanning more than one registry.

\section{Results} \label{sec: Results}

\subsection{Comparison with varying arrival rates in registries}

An important consideration in multi-registry kidney exchange programs (mKEP) is understanding how the size of participating registries influences the potential benefits of collaboration. Specifically, the gains from participation may vary depending on registry size, with smaller registries potentially benefiting more in size-imbalanced scenarios. Estimating the benefits across different registry sizes is therefore crucial for informing decisions about joining mKEP initiatives. In this analysis, we examine a two-registry case to compare the outcomes for registries of varying sizes. We present results on the number of patients matched and provide a detailed comparison of blood group distributions to assess the nature and extent of the benefits realized.

\begin{table}[ht]
\centering
\footnotesize
\renewcommand{\arraystretch}{1.15}
\setlength{\tabcolsep}{3pt}
\resizebox{\textwidth}{!}{%
\begin{tabular}{l l cccc cccc}
\toprule
\multirow{2}{*}{Scenario} & \multirow{2}{*}{Model} & \multicolumn{4}{c}{Registry 1} & \multicolumn{4}{c}{Registry 2} \\
\cmidrule(lr){3-6} \cmidrule(lr){7-10}
 & & Tx (SD) & ES (SD) & $\Delta$Tx & $\Delta$ES & Tx (SD) & ES (SD) & $\Delta$Tx & $\Delta$ES \\
\midrule
\multirow{4}{*}{\makecell[l]{Baseline\\ $R_1,R_2\sim\mathrm{Pois}(5)$}}
 & Individual  & 49.5 (9.1) & 56.3 (2.5) & --      & --     & 49.3 (10.6) & 56.3 (2.6) & --      & --     \\
 & mKEP        & 54.5 (7.1) & 56.0 (2.8) & +10.3\% & -0.5\% & 54.0 (9.3)  & 55.9 (2.3) & +9.5\%  & -0.7\% \\
 & mKEP--IR    & 54.0 (7.0) & 55.9 (2.4) & +9.1\%  & -0.7\% & 54.4 (9.9)  & 56.1 (2.4) & +10.3\% & -0.4\% \\
 & mKEP--Fair  & 54.0 (8.2) & 56.1 (2.6) & +9.2\%  & -0.3\% & 54.1 (10.2) & 56.1 (2.5) & +9.7\%  & -0.4\% \\
\midrule
\multirow{4}{*}{\makecell[l]{Arrival asymmetry\\ $R_1\sim\mathrm{Pois}(5)$,\\ $R_2\sim\mathrm{Pois}(10)$}}
 & Individual  & 50.1 (9.2) & 56.2 (2.4) & --      & --     & 106.8 (12.5) & 55.9 (1.5) & --     & --     \\
 & mKEP        & 56.5 (8.6) & 56.4 (2.1) & +12.8\% & +0.3\% & 110.5 (11.0) & 56.3 (1.6) & +3.5\% & +0.7\% \\
 & mKEP--IR    & 56.1 (8.8) & 55.9 (2.3) & +12.0\% & -0.5\% & 110.1 (12.1) & 56.5 (1.7) & +3.2\% & +1.1\% \\
 & mKEP--Fair  & 54.2 (9.3) & 56.4 (2.1) & +8.2\%  & +0.2\% & 111.6 (11.9) & 56.2 (1.6) & +4.5\% & +0.5\% \\
\midrule
\multirow{4}{*}{\makecell[l]{Both high\\ $R_1,R_2\sim\mathrm{Pois}(10)$}}
 & Individual  & 108.4 (13.2) & 56.3 (1.5) & --     & --     & 107.8 (11.7) & 55.9 (2.0) & --     & --     \\
 & mKEP        & 115.4 (12.2) & 55.9 (1.4) & +6.5\% & -0.7\% & 113.4 (10.6) & 56.4 (1.5) & +5.2\% & +0.9\% \\
 & mKEP--IR    & 114.8 (13.1) & 56.1 (1.6) & +5.9\% & -0.5\% & 113.9 (10.5) & 56.3 (1.6) & +5.7\% & +0.7\% \\
 & mKEP--Fair  & 114.3 (13.0) & 55.8 (1.5) & +5.4\% & -0.9\% & 113.3 (11.4) & 55.9 (1.6) & +5.1\% & +0.0\% \\
\bottomrule
\end{tabular}%
}
\caption{Effect of arrival rate -- both its level and its asymmetry across
registries -- on transplant volume (Tx) and match quality (ES), by registry
and allocation model. Standard deviations (in parentheses) are computed across
$N=50$ replications of 24 rounds each. $\Delta$Tx and $\Delta$ES report each
pooled model's percentage change relative to that registry's own Individual
value in the same scenario; positive values indicate improvement over acting
alone. All other parameters held at baseline (ASTRA blood-group distribution
for both registries, $b_1=b_2=3$, $B=5$, dropout probability 0.2 for both).}
\label{tab:results_arrival}
\end{table}
 
\begin{table}[ht]
\centering
\footnotesize
\renewcommand{\arraystretch}{1.15}
\setlength{\tabcolsep}{3pt}
\resizebox{\textwidth}{!}{%
\begin{tabular}{l l cccc cccc}
\toprule
\multirow{2}{*}{Scenario} & \multirow{2}{*}{Model} & \multicolumn{4}{c}{Registry 1} & \multicolumn{4}{c}{Registry 2} \\
\cmidrule(lr){3-6} \cmidrule(lr){7-10}
 & & Dropouts (SD) & Wait & $\Delta$Drop & $\Delta$Wait & Dropouts (SD) & Wait & $\Delta$Drop & $\Delta$Wait \\
\midrule
\multirow{4}{*}{\makecell[l]{Baseline\\ $R_1,R_2\sim\mathrm{Pois}(5)$}}
 & Individual  & 62.5 (9.1) & 1.68 & --     & --      & 60.9 (8.8) & 1.68 & --     & --      \\
 & mKEP        & 58.0 (7.6) & 1.57 & -7.2\% & -6.8\%  & 56.5 (7.4) & 1.50 & -7.2\% & -10.6\% \\
 & mKEP--IR    & 58.7 (8.4) & 1.52 & -6.0\% & -9.6\%  & 56.0 (7.6) & 1.49 & -8.0\% & -11.5\% \\
 & mKEP--Fair  & 58.5 (8.6) & 1.52 & -6.4\% & -9.8\%  & 56.5 (7.0) & 1.48 & -7.3\% & -11.9\% \\
\midrule
\multirow{4}{*}{\makecell[l]{Arrival asymmetry\\ $R_1\sim\mathrm{Pois}(5)$,\\ $R_2\sim\mathrm{Pois}(10)$}}
 & Individual  & 60.5 (8.4) & 1.69 & --      & --      & 114.9 (14.0) & 1.56 & --     & --      \\
 & mKEP        & 55.0 (8.3) & 1.45 & -9.1\%  & -13.9\% & 112.0 (13.2) & 1.47 & -2.5\% & -5.8\%  \\
 & mKEP--IR    & 55.1 (8.3) & 1.46 & -8.8\%  & -13.3\% & 112.5 (14.2) & 1.46 & -2.1\% & -6.5\%  \\
 & mKEP--Fair  & 56.8 (8.3) & 1.48 & -6.1\%  & -12.6\% & 111.2 (13.3) & 1.45 & -3.2\% & -7.0\%  \\
\midrule
\multirow{4}{*}{\makecell[l]{Both high\\ $R_1,R_2\sim\mathrm{Pois}(10)$}}
 & Individual  & 116.4 (10.5) & 1.53 & --     & --      & 117.0 (12.1) & 1.56 & --     & --      \\
 & mKEP        & 111.5 (9.8)  & 1.42 & -4.3\% & -7.3\%  & 112.3 (11.7) & 1.39 & -4.1\% & -10.7\% \\
 & mKEP--IR    & 111.8 (9.7)  & 1.39 & -4.0\% & -8.8\%  & 111.4 (13.2) & 1.41 & -4.8\% & -9.7\%  \\
 & mKEP--Fair  & 112.1 (10.1) & 1.38 & -3.7\% & -9.4\%  & 112.4 (13.0) & 1.40 & -4.0\% & -10.1\% \\
\bottomrule
\end{tabular}%
}
\caption{Total dropouts over the 24-round horizon (mean, SD across $N=50$
replications in parentheses) and average waiting time in rounds until match,
for the same three arrival-rate scenarios and four allocation models as
Table~\ref{tab:results_arrival}. Waiting time is computed only over matched
patients; it does not include patients who dropped out unmatched. $\Delta$Drop
and $\Delta$Wait report each pooled model's percentage change relative to
Individual in the same scenario; unlike Table~\ref{tab:results_arrival},
\emph{negative} values here indicate improvement (fewer dropouts, shorter
waits) over acting alone.}
\label{tab:results_arrival_dropwait}
\end{table}

Table~\ref{tab:results_arrival} compares multi-registry kidney exchange (mKEP) allocation against individually operating registries under varying arrival rates. At baseline, when both registries have low and equal arrival rates, pooling improves transplant volume substantially and similarly across all three pooled mechanisms (roughly 9--10\%) for both registries, while match quality is essentially unchanged: the edge-score (ES) differences are small (within about 1\%) and well inside one standard deviation, so they should be read as no detectable change in quality rather than a quality loss. When one registry has a markedly larger arrival rate than the other, the relative benefit of pooling becomes more pronounced for the lower-arrival registry, whose transplant volume rises by 8--13\% depending on mechanism, with a small, mechanism-dependent change in match quality. The higher-arrival registry instead sees a much smaller transplant-volume gain (roughly 3--5\%) but a more consistently positive, if still modest, change in match quality across all three mechanisms. When both registries have high arrival rates, the transplant-volume benefit of pooling moderates to roughly 5--6\% for both registries, while match quality remains close to the individual baseline for both.

Comparing the three pooled mechanisms directly under arrival asymmetry, the Shapley-fair allocation is the only one that narrows the gap in relative benefit between the two registries: the difference in transplant-volume gain between the lower- and higher-arrival registry is 9.3 and 8.8 percentage points under the unconstrained pool and mKEP--IR respectively, but only 3.7 percentage points under mKEP--Fair. This narrowing, however, is achieved by redistributing gain toward the higher-arrival registry rather than by both registries doing better: mKEP--IR gives the lower-arrival registry a larger transplant-volume gain than mKEP--Fair does (+12.0\% vs.\ +8.2\%), while mKEP--Fair gives the higher-arrival registry a larger gain than mKEP--IR does (+4.5\% vs.\ +3.2\%). Which mechanism is preferable for a given registry therefore depends on which side of the arrival-rate asymmetry that registry sits on, rather than one mechanism uniformly dominating the other.

Table~\ref{tab:results_arrival_dropwait} reports the corresponding effect of arrival rate on dropouts and waiting time, and shows the same qualitative pattern as transplant volume. At baseline, pooling reduces both dropouts (by roughly 6--8\%) and waiting time (by roughly 7--12\%) across all three mechanisms, with the reduction in waiting time consistently exceeding the reduction in dropouts. Under arrival asymmetry, this benefit becomes substantially larger for the lower-arrival registry, whose dropouts and waiting time both fall by 6--14\%, while the higher-arrival registry sees a much smaller reduction in dropouts (2--3\%) alongside a still-meaningful reduction in waiting time (6--7\%). When both registries have high arrival rates, the reduction in dropouts and waiting time settles to an intermediate
level (roughly 4--5\% and 7--11\%, respectively) for both registries. Directly comparing mechanisms under arrival asymmetry reproduces the same redistribution pattern found for transplant volume: mKEP--IR achieves a larger dropout reduction for the lower-arrival registry (-8.8\% vs.\ -6.1\% under mKEP--Fair), while mKEP--Fair achieves a larger reduction for the higher-arrival registry (-3.2\% vs.\ -2.1\% under mKEP--IR). This is consistent with the transplant-volume result and reinforces that the fairness mechanism systematically redistributes benefit toward whichever registry gains least under individual rationality or unconstrained pooling, across transplant volume, dropouts, and waiting time alike. Standard deviations are not currently reported for waiting time, so the smaller waiting-time differences between mechanisms in particular should be treated
as suggestive rather than statistically confirmed.

\begin{figure}[p]
\centering
\begin{subfigure}[t]{0.8\textwidth}
    \centering
    \includegraphics[width=\textwidth]{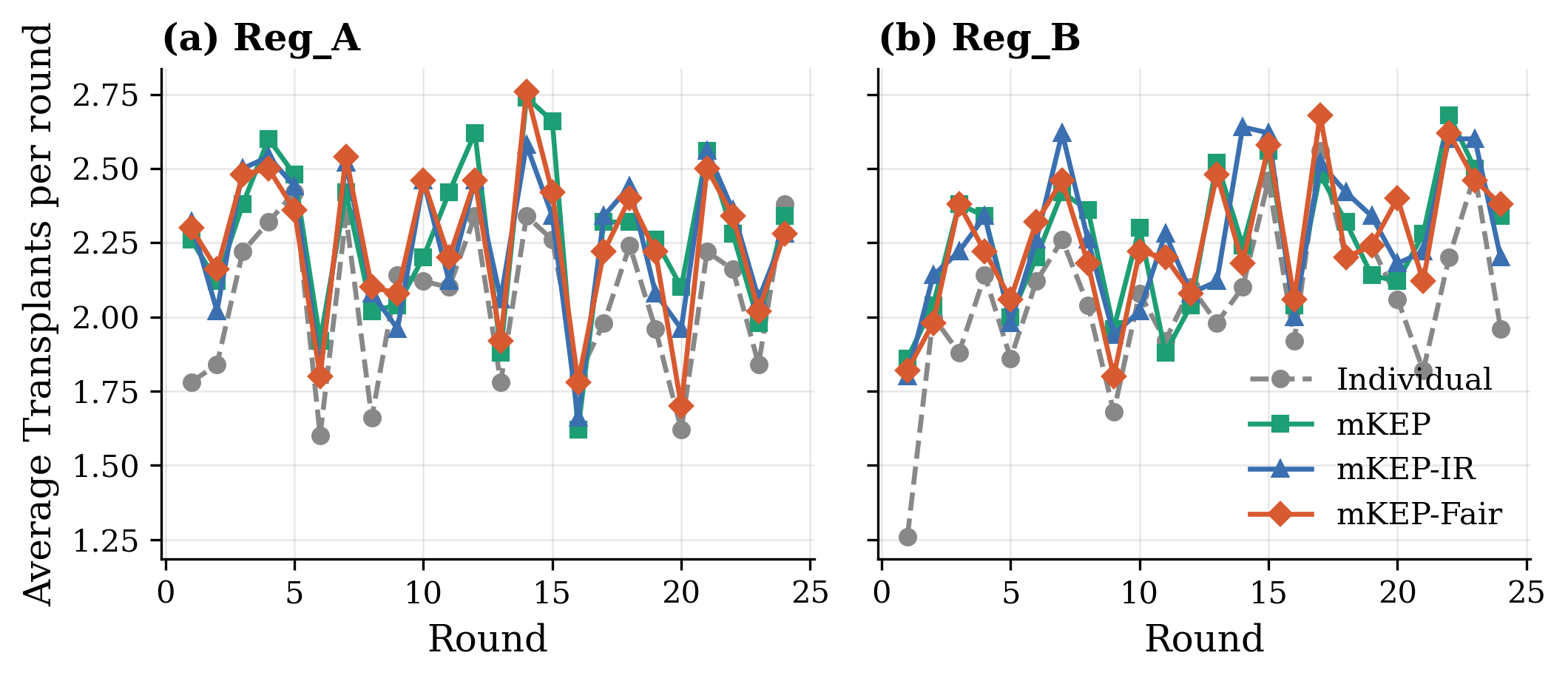}
    \caption{Baseline ($R_1,R_2\sim\mathrm{Pois}(5)$)}
    \label{fig:tx_arrival_baseline}
\end{subfigure}
\begin{subfigure}[t]{0.8\textwidth}
    \centering
    \includegraphics[width=\textwidth]{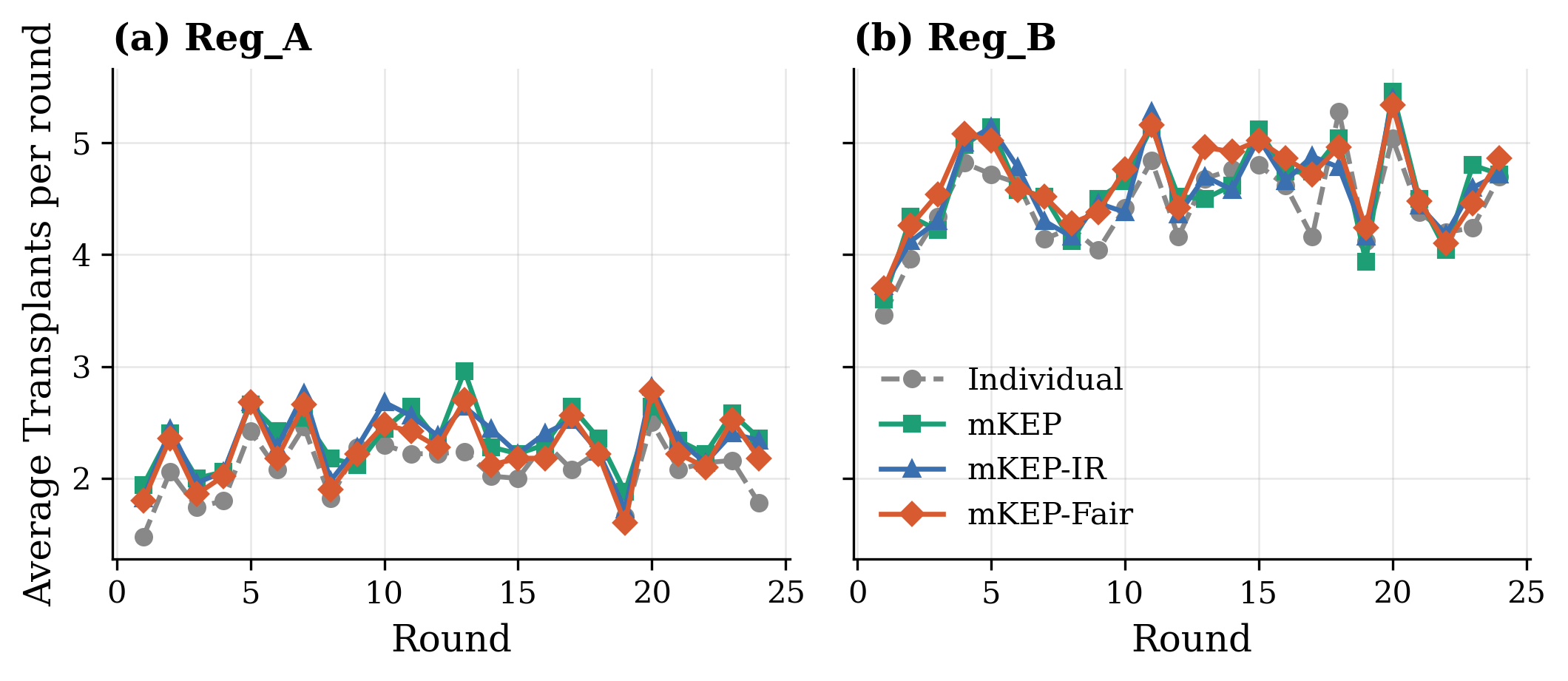}
    \caption{Arrival asymmetry ($R_1\sim\mathrm{Pois}(5)$, $R_2\sim\mathrm{Pois}(10)$)}
    \label{fig:tx_arrival_asym}
\end{subfigure}
\begin{subfigure}[t]{0.8\textwidth}
    \centering
    \includegraphics[width=\textwidth]{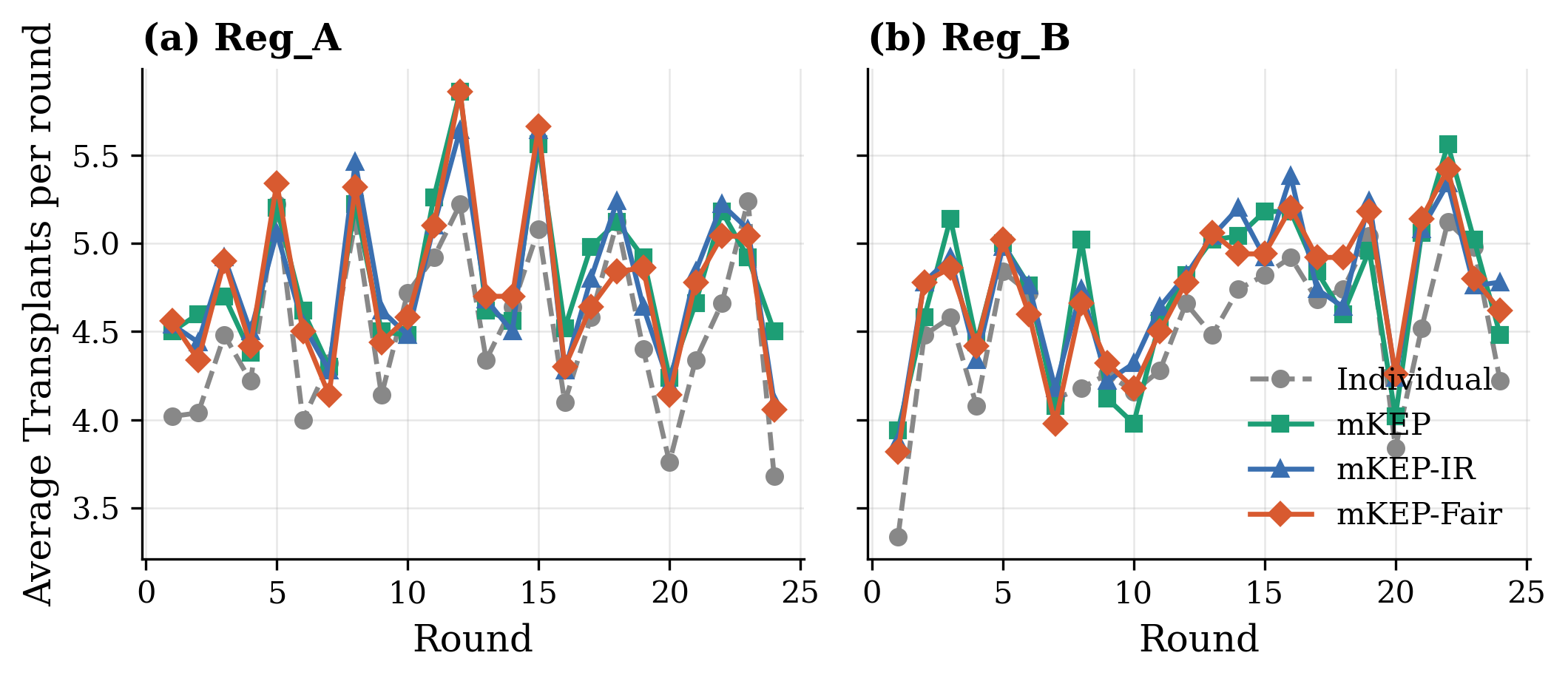}
    \caption{Both high ($R_1,R_2\sim\mathrm{Pois}(10)$)}
    \label{fig:tx_arrival_bothhigh}
\end{subfigure}
\caption{Average transplants per round over the 24-round horizon, by registry
and allocation model, across the three arrival-rate scenarios of
Table~\ref{tab:results_arrival}. Dashed gray line shows the Individual
(standalone) baseline for reference against the three pooled mechanisms.}
\label{fig:tx_arrival_all}
\end{figure}

\begin{figure}[p]
\centering
\begin{subfigure}[t]{0.8\textwidth}
    \centering
    \includegraphics[width=\textwidth]{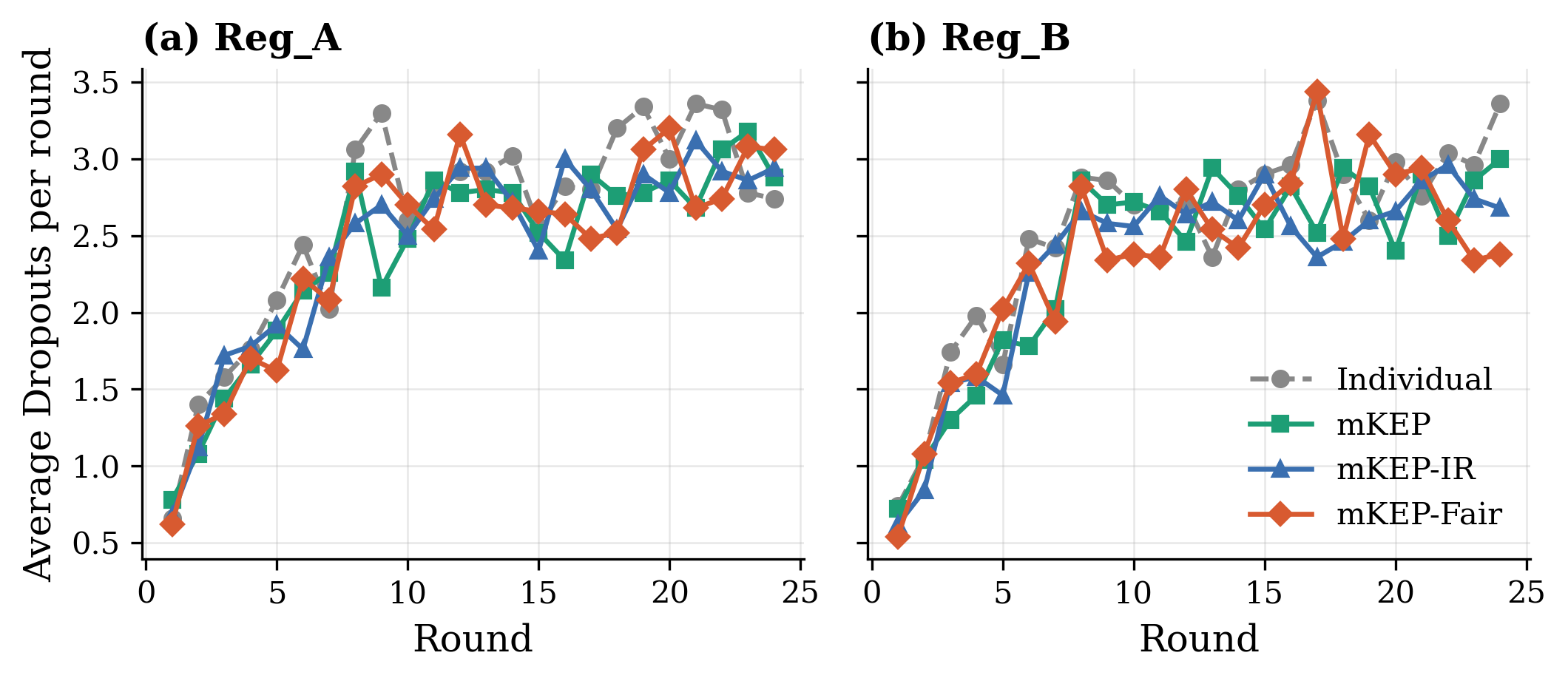}
    \caption{Baseline ($R_1,R_2\sim\mathrm{Pois}(5)$)}
    \label{fig:drop_arrival_baseline}
\end{subfigure}
\vspace{0.5em}
\begin{subfigure}[t]{0.8\textwidth}
    \centering
    \includegraphics[width=\textwidth]{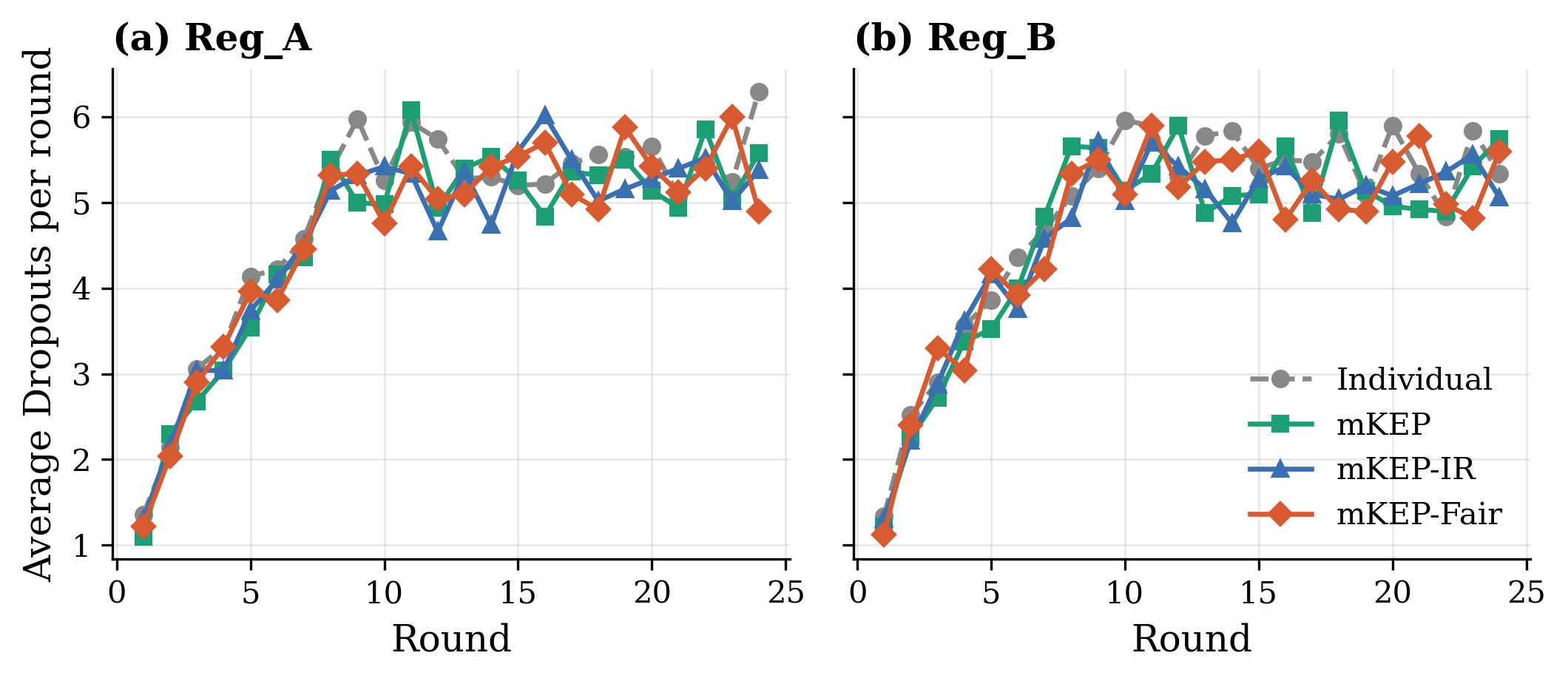}
    \caption{Arrival asymmetry ($R_1\sim\mathrm{Pois}(5)$, $R_2\sim\mathrm{Pois}(10)$)}
    \label{fig:drop_arrival_asym}
\end{subfigure}
\vspace{0.5em}
\begin{subfigure}[t]{0.8\textwidth}
    \centering
    \includegraphics[width=\textwidth]{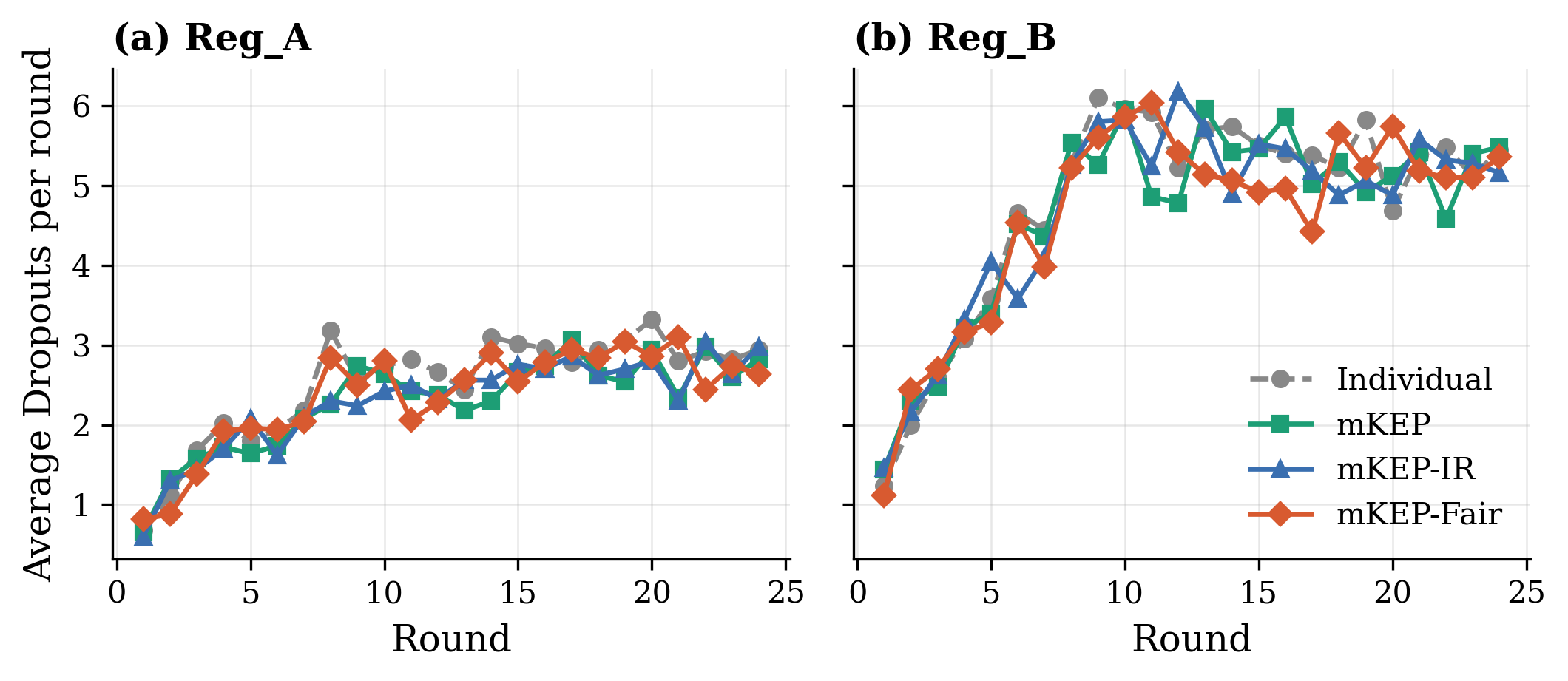}
    \caption{Both high ($R_1,R_2\sim\mathrm{Pois}(10)$)}
    \label{fig:drop_arrival_bothhigh}
\end{subfigure}
\caption{Average dropouts per round over the 24-round horizon, by registry and
allocation model, across the three arrival-rate scenarios of
Table~\ref{tab:results_arrival}.}
\label{fig:drop_arrival_all}
\end{figure}

Figure~\ref{fig:tx_arrival_all} shows the average number of transplants per round over the full 24-round horizon for each of the three arrival-rate scenarios. Two patterns emerge that are not visible from the summary statistics in Table~\ref{tab:results_arrival} alone. First, in the two symmetric scenarios (baseline and both-high), the three pooled mechanisms track closely together and remain above the Individual curve in nearly every round, not just on average -- indicating that the transplant-volume benefit of pooling is a consistent, round-by-round effect rather than one driven by a
handful of favorable rounds. Second, under arrival asymmetry, this separation becomes markedly less consistent for the higher-arrival registry: its Individual curve intermittently matches or exceeds the pooled curves (for example, around round 20), whereas the lower-arrival registry's pooled curves remain clearly and consistently above Individual throughout. This suggests that arrival-rate asymmetry does not only shrink the higher-arrival registry's average benefit from pooling, as already seen in Table~\ref{tab:results_arrival}, but also makes that benefit less reliable from round to round.

Figure~\ref{fig:drop_arrival_all} shows the corresponding dropout
trajectories. All curves rise over the first 8--10 rounds before leveling off, reflecting the pool filling from empty at the start of the simulation and dropouts reaching a steady state once the patient population stabilizes; results from these early rounds should accordingly be interpreted with more caution than the steady-state portion of the horizon. Once stabilized, the ordering mirrors the transplant-volume figure: Individual sits above the
three pooled mechanisms in the symmetric scenarios, with the gap narrowing and becoming less consistent under arrival asymmetry, particularly for the higher-arrival registry.

\subsection{Comparison with varying blood group distribution}

A second key dimension in evaluating the effectiveness of multi-registry kidney exchange
programs (mKEP) is understanding how variations in blood group distributions across registries
affect matching outcomes. In practice, the composition of donor and recipient blood types within
a registry can vary significantly depending on regional demographics, healthcare access, and
participation patterns. Such heterogeneity can have important implications for both the
feasibility and equity of matches in an mKEP framework. To investigate this, we simulate matching
outcomes using two empirically observed blood type distributions: one derived from the ASTRA
registry and the other from the APKD registry.

For kidney exchange programs to operate optimally, there should ideally be a close alignment
between donor and recipient blood group distributions. This balance minimizes the risk of having
unmatchable patients due to shortages of specific blood types. However, real-world kidney exchange
registries often reflect imbalances driven by biological and behavioral factors. For instance,
O-type patients who can only receive kidneys from O-type donors tend to be overrepresented
on waitlists, whereas AB-type donors who can donate to any recipient are often more
prevalent in donor pools. These patterns result from the natural constraints of blood group
compatibility and the selective enrollment of incompatible pairs.

The APKD registry, in particular, contains a significantly higher proportion of O-type donors
compared to the ASTRA registry: O-type donors make up 34.5\% of the APKD donor pool, against only
3.3\% of the ASTRA donor pool, a roughly tenfold difference, giving APKD a higher proportion of
easy-to-match pairs than ASTRA. This distinction provides a valuable
opportunity to examine how such differences influence the overall performance and equity of mKEP. The specific distributions employed in this analysis are presented in
\ref{sec: Data_distribution}, offering a detailed reference for understanding the
simulated scenarios.

Table~\ref{tab:results_bloodgroup} reports transplant volume and match quality under this
composition asymmetry, crossed against four arrival-rate combinations, with $R_1$ fixed to the
ASTRA distribution and $R_2$ to the APKD distribution throughout. Registry 1 gains substantially
under every pooled mechanism across all four scenarios, with transplant volume rising by
9--45\% depending on mechanism and arrival combination. Registry 2's transplant volume, in
contrast, \emph{falls} under the unconstrained pool in all four scenarios (by 0.4--5.5\%) --
the only registry, across every scenario examined in this study, for which unconstrained pooling
produces a net loss relative to acting alone. Both mKEP--IR and mKEP--Fair reverse this: Registry
2's transplant volume is positive under both mechanisms in every scenario. When Registry 1
(ASTRA) has the higher arrival rate and Registry 2 (APKD) the lower one, mKEP--IR provides a
better outcome for Registry 2 than mKEP--Fair does (+6.2\% vs.\ +5.2\%). Match quality changes
are modest throughout, mostly within about one percentage point in either direction, for both
registries and all three pooled mechanisms.
 
The size of Registry 2's gain under mKEP--IR and mKEP--Fair, however, depends strongly on whether
Registry 2 is also the higher-arrival registry in a given scenario, and not on its blood-group
composition alone. When Registry 2 has the lower or equal arrival rate, its gain under IR or Fair
ranges from 4.6\% to 6.2\%; when Registry 2 instead has the higher arrival rate, that gain falls
to 0.6--2.6\%. In absolute terms, this smaller gain amounts to only 2--4 additional transplants
over the full six-year horizon for a registry already realizing roughly 159--167 transplants on
its own. For instance, when both registries have the higher arrival rate, Registry 2 rises from
166.6 (Individual) to 168.9 under mKEP--IR and 170.9 under mKEP--Fair. Registry 2's standard
deviation in transplant volume is approximately 17--18 in these scenarios, so a shift of this
size is well within one standard deviation.

Taken together, these results suggest that a registry that is both easy to match internally and already high-volume has comparatively little to gain from joining an mKEP, even under a fairness mechanism explicitly designed to protect it, and may reasonably be hesitant to bear the coordination and data-sharing costs of participation for a benefit this small whereas a smaller or harder-to-match registry,
such as Registry 1 throughout these scenarios, stands to gain substantially more.
 
The blood-type breakdown in Table~\ref{tab:results_bloodgroup_bytype} shows that Registry 2's
aggregate loss under unconstrained pooling is driven almost entirely by O-type recipients. Its
A-, B-, and AB-type matches rise under every mechanism and every scenario, while its O-type
matches fall by 9--15\% under the unconstrained pool. Registry 1's O-type matches move sharply in
the opposite direction (+98\% to +357\%), consistent with pooling reallocating Registry 2's
comparatively large O-type donor supply toward Registry 1's O-type recipients, who have little
access to O-type donors internally. Both mKEP--IR and mKEP--Fair substantially narrow this loss
for Registry 2, but do not always eliminate it: Registry 2's O-type matches remain below their
standalone level in three of the four scenarios under IR (by up to 2.6\%) and in three of the
four under Fair (by up to 1.0\%), with only one scenario under each mechanism yielding a small
positive change. A, B, and AB recipients, by contrast, benefit under every pooled mechanism in
every scenario, for both registries.

\begin{table}[ht]
\centering
\footnotesize
\renewcommand{\arraystretch}{1.15}
\setlength{\tabcolsep}{3pt}
\resizebox{\textwidth}{!}{%
\begin{tabular}{l l cccc cccc}
\toprule
\multirow{2}{*}{Scenario} & \multirow{2}{*}{Model} & \multicolumn{4}{c}{Registry 1} & \multicolumn{4}{c}{Registry 2} \\
\cmidrule(lr){3-6} \cmidrule(lr){7-10}
 & & Tx (SD) & ES (SD) & $\Delta$Tx & $\Delta$ES & Tx (SD) & ES (SD) & $\Delta$Tx & $\Delta$ES \\
\midrule
\multirow{4}{*}{\makecell[l]{BG asymmetry, equal arrivals\\ $R_1,R_2\sim\mathrm{Pois}(5)$\\ $R_1$: ASTRA, $R_2$: APKD}}
 & Individual  & 48.4 (9.1) & 56.1 (2.7) & --      & --     & 75.0 (10.4) & 56.4 (1.8) & --      & --     \\
 & mKEP        & 60.8 (8.7) & 56.6 (2.4) & +25.5\% & +0.9\% & 74.7 (8.6)  & 56.8 (2.3) & -0.4\%  & +0.6\% \\
 & mKEP--IR    & 56.9 (8.7) & 56.5 (2.4) & +17.5\% & +0.7\% & 78.4 (9.7)  & 56.3 (2.4) & +4.6\%  & -0.2\% \\
 & mKEP--Fair  & 56.0 (9.0) & 56.4 (2.4) & +15.7\% & +0.5\% & 78.9 (10.4) & 56.2 (1.9) & +5.2\%  & -0.5\% \\
\midrule
\multirow{4}{*}{\makecell[l]{BG asymmetry, unequal arrivals\\ $R_1\sim\mathrm{Pois}(5)$, $R_2\sim\mathrm{Pois}(10)$\\ $R_1$: ASTRA, $R_2$: APKD}}
 & Individual  & 47.7 (7.0) & 56.4 (2.5) & --      & --     & 159.1 (18.1) & 56.1 (1.4) & --     & --     \\
 & mKEP        & 69.2 (7.8) & 56.1 (2.4) & +45.1\% & -0.5\% & 152.3 (17.0) & 56.1 (1.5) & -4.3\% & -0.1\% \\
 & mKEP--IR    & 61.0 (7.0) & 56.4 (2.2) & +27.9\% & +0.0\% & 160.1 (17.7) & 56.0 (1.3) & +0.6\% & -0.2\% \\
 & mKEP--Fair  & 57.1 (7.0) & 56.7 (2.2) & +19.8\% & +0.4\% & 162.0 (18.5) & 56.2 (1.5) & +1.8\% & +0.1\% \\
\midrule
\multirow{4}{*}{\makecell[l]{BG asymmetry, unequal arrivals (swapped)\\ $R_1\sim\mathrm{Pois}(10)$, $R_2\sim\mathrm{Pois}(5)$\\ $R_1$: ASTRA, $R_2$: APKD}}
 & Individual  & 106.0 (13.0) & 55.8 (1.8) & --     & --     & 77.3 (10.4) & 56.2 (1.8) & --     & --     \\
 & mKEP        & 120.5 (12.5) & 55.8 (1.9) & +13.7\% & +0.0\% & 76.6 (8.3)  & 55.8 (1.9) & -0.9\% & -0.7\% \\
 & mKEP--IR    & 115.2 (12.4) & 56.0 (1.9) & +8.7\%  & +0.4\% & 82.1 (9.5)  & 56.2 (2.2) & +6.2\% & +0.0\% \\
 & mKEP--Fair  & 115.7 (12.4) & 56.0 (1.8) & +9.1\%  & +0.4\% & 81.3 (9.6)  & 55.9 (2.1) & +5.2\% & -0.5\% \\
\midrule
\multirow{4}{*}{\makecell[l]{BG asymmetry, equal (higher) arrivals\\ $R_1,R_2\sim\mathrm{Pois}(10)$\\ $R_1$: ASTRA, $R_2$: APKD}}
 & Individual  & 108.0 (13.0) & 56.4 (1.8) & --      & --     & 166.6 (17.2) & 56.1 (1.2) & --     & --     \\
 & mKEP        & 134.4 (13.4) & 56.2 (1.7) & +24.5\% & -0.3\% & 157.5 (15.5) & 56.2 (1.4) & -5.5\% & +0.2\% \\
 & mKEP--IR    & 122.6 (12.8) & 56.3 (1.5) & +13.5\% & -0.2\% & 168.9 (17.0) & 56.4 (1.2) & +1.3\% & +0.6\% \\
 & mKEP--Fair  & 119.1 (12.8) & 56.4 (1.4) & +10.3\% & +0.2\% & 170.9 (17.3) & 56.3 (1.5) & +2.6\% & +0.5\% \\
\bottomrule
\end{tabular}%
}
\caption{Effect of blood-group-composition asymmetry ($R_1$: ASTRA,
$R_2$: APKD, fixed throughout) on transplant volume (Tx) and match quality
(ES), crossed against four arrival-rate combinations, by registry and
allocation model. The second and third scenario blocks form a matched pair
that swaps \emph{which} registry (and hence which blood-group distribution)
carries the higher arrival rate. $\Delta$Tx and $\Delta$ES report each
pooled model's percentage change relative to that registry's own Individual
value in the same scenario, as in Table~\ref{tab:results_arrival}. Standard
deviations (in parentheses) are computed across $N=50$ replications of 24
rounds each. All other parameters held at baseline ($b_1=b_2=3$, $B=5$,
dropout probability 0.2 for both).}
\label{tab:results_bloodgroup}
\end{table}
 
\begin{table}[ht]
\centering
\footnotesize
\renewcommand{\arraystretch}{1.15}
\setlength{\tabcolsep}{2.5pt}
\resizebox{\textwidth}{!}{%
\begin{tabular}{l l cccc cccc}
\toprule
\multirow{2}{*}{Scenario} & \multirow{2}{*}{Model} & \multicolumn{4}{c}{Registry 1 (ASTRA)} & \multicolumn{4}{c}{Registry 2 (APKD)} \\
\cmidrule(lr){3-6} \cmidrule(lr){7-10}
 & & A & B & AB & O & A & B & AB & O \\
\midrule
\multirow{4}{*}{\makecell[l]{Equal arrivals\\ $R_1,R_2\sim\mathrm{Pois}(5)$}}
 & Individual  & 21.44 & 21.70 & 2.00 & 3.28 & 23.70 & 11.58 & 2.10 & 37.60 \\
 & mKEP        & 24.04 (+12\%) & 24.42 (+13\%) & 2.22 (+11\%) & 10.10 (+207.9\%) & 24.98 (+5\%) & 13.32 (+15\%) & 2.24 (+7\%) & 34.12 (-9.3\%) \\
 & mKEP--IR    & 24.10 (+12\%) & 23.38 (+8\%)  & 2.22 (+11\%) & 7.20 (+119.5\%)  & 25.18 (+6\%) & 13.92 (+20\%) & 2.24 (+7\%) & 37.10 (-1.3\%) \\
 & mKEP--Fair  & 23.86 (+11\%) & 23.28 (+7\%)  & 2.20 (+10\%) & 6.68 (+103.7\%)  & 25.46 (+7\%) & 13.76 (+19\%) & 2.24 (+7\%) & 37.44 (-0.4\%) \\
\midrule
\multirow{4}{*}{\makecell[l]{Unequal arrivals\\ $R_1\sim\mathrm{Pois}(5)$,\\ $R_2\sim\mathrm{Pois}(10)$}}
 & Individual  & 21.08 & 21.34 & 1.88 & 3.38 & 51.78 & 25.96 & 4.78 & 76.60 \\
 & mKEP        & 25.30 (+20\%) & 26.28 (+23\%) & 2.16 (+15\%) & 15.44 (+356.8\%) & 52.16 (+1\%) & 27.20 (+5\%) & 4.88 (+2\%) & 68.06 (-11.1\%) \\
 & mKEP--IR    & 24.66 (+17\%) & 25.34 (+19\%) & 2.18 (+16\%) & 8.80 (+160.4\%)  & 52.76 (+2\%) & 27.84 (+7\%) & 4.94 (+3\%) & 74.60 (-2.6\%)  \\
 & mKEP--Fair  & 24.06 (+14\%) & 24.24 (+14\%) & 2.10 (+12\%) & 6.70 (+98.2\%)   & 52.78 (+2\%) & 27.88 (+7\%) & 4.86 (+2\%) & 76.50 (-0.1\%)  \\
\midrule
\multirow{4}{*}{\makecell[l]{Unequal arrivals (swapped)\\ $R_1\sim\mathrm{Pois}(10)$,\\ $R_2\sim\mathrm{Pois}(5)$}}
 & Individual  & 47.14 & 47.40 & 4.04 & 7.44 & 24.36 & 12.44 & 2.14 & 38.36 \\
 & mKEP        & 49.90 (+6\%) & 50.54 (+7\%) & 4.24 (+5\%) & 15.84 (+112.9\%) & 26.10 (+7\%) & 14.42 (+16\%) & 2.30 (+8\%) & 33.82 (-11.8\%) \\
 & mKEP--IR    & 50.08 (+6\%) & 49.56 (+5\%) & 4.22 (+5\%) & 11.38 (+53.0\%)  & 26.44 (+9\%) & 14.70 (+18\%) & 2.30 (+8\%) & 38.62 (+0.7\%)  \\
 & mKEP--Fair  & 49.88 (+6\%) & 49.76 (+5\%) & 4.20 (+4\%) & 11.82 (+58.9\%)  & 26.34 (+8\%) & 14.70 (+18\%) & 2.30 (+8\%) & 37.98 (-1.0\%)  \\
\midrule
\multirow{4}{*}{\makecell[l]{Equal (higher) arrivals\\ $R_1,R_2\sim\mathrm{Pois}(10)$}}
 & Individual  & 48.24 & 48.18 & 3.56 & 8.02 & 52.52 & 27.94 & 5.14 & 81.02 \\
 & mKEP        & 52.64 (+9\%)  & 53.14 (+10\%) & 3.76 (+6\%) & 24.90 (+210.5\%) & 53.48 (+2\%) & 29.66 (+6\%) & 5.26 (+2\%) & 69.10 (-14.7\%) \\
 & mKEP--IR    & 52.96 (+10\%) & 51.32 (+7\%)  & 3.82 (+7\%) & 14.50 (+80.8\%)  & 53.78 (+2\%) & 30.34 (+9\%) & 5.22 (+2\%) & 79.52 (-1.9\%)  \\
 & mKEP--Fair  & 52.38 (+9\%)  & 50.36 (+5\%)  & 3.78 (+6\%) & 12.56 (+56.6\%)  & 53.92 (+3\%) & 30.38 (+9\%) & 5.28 (+3\%) & 81.34 (+0.4\%)  \\
\bottomrule
\end{tabular}%
}
\caption{Average number of matched pairs per replication by recipient blood
type, registry, and allocation model, across the same four blood-group-asymmetry scenarios as Table~\ref{tab:results_bloodgroup} ($R_1$: ASTRA,
$R_2$: APKD). Percentage in parentheses is change relative to that registry's
own Individual value in the same scenario, as elsewhere in this study. AB
counts are small in absolute terms (2--5 matches), so the associated
percentage changes should be interpreted with additional caution.}
\label{tab:results_bloodgroup_bytype}
\end{table}

\subsection{Comparison with varying bound on cycle lengths}

The third factor that may influence mKEP outcomes is the constraint on cycle length imposed by each
participating registry. These cycle length bounds can vary due to several operational and
logistical factors, such as whether the registry is managed centrally or across multiple centers,
the transplant capacity of participating hospitals, and whether transplants are performed
simultaneously or non-simultaneously. Here, we compare the performance of mKEP to that of
individual registry solutions under different cycle length constraints. The analysis is based on
blood group distributions from the ASTRA registry, with lower arrival rates used to reflect more
conservative operational conditions.

Table~\ref{tab:results_bound} reports outcomes under an asymmetric, tighter cycle-length bound
for Registry 2 ($b_1=3$, $b_2=2$) together with a tightened cross-registry bound ($B=3$), crossed
against three arrival-rate combinations. 
Relative to the equivalent low-arrival scenario under the relaxed, symmetric bound in
Table~\ref{tab:results_arrival} ($b_1=b_2=3$, $B=5$), tightening the bound increases the relative
benefit of pooling for both registries under all three mechanisms. At equal arrival rates, Registry
1's gain rises from 9.1--9.2\% under the relaxed bound to 12.1--12.3\% under the tighter bound,
and Registry 2's from 9.7--10.3\% to 12.1--13.3\%. The effect is more pronounced under arrival
asymmetry: Registry 1's gain rises from 8.2--12.0\% under the relaxed bound to 19.4--21.8\% under
the tighter bound. This is consistent with a tighter cycle-length bound constraining what a
registry can achieve on its own more than it constrains what the same registry can achieve once
pooled: even a tightened cross-registry bound typically still admits more matching configurations
across a larger, combined graph than the same bound permits within a single smaller registry, so
pooling recovers a larger share of the compatibility that a tight individual bound otherwise
leaves unrealized.

Comparing mKEP--IR and mKEP--Fair directly, mKEP--IR gives a larger gain than mKEP--Fair to
Registry 2, the registry operating under the tighter internal bound in all three scenarios,
though the margin varies considerably (0.2 to 1.4 percentage points) and narrows sharply when
Registry 2 also has the higher arrival rate. This does not reproduce the clean crossover seen
under arrival asymmetry alone in Table~\ref{tab:results_arrival}, where mKEP--IR and mKEP--Fair
favored opposite registries depending on which had the higher arrival rate: here, mKEP--IR gives
the larger gain to Registry 2 in every scenario, including the one in which Registry 1 has the
lower arrival rate. 
Match-quality changes remain modest throughout, generally within about one percentage point in
either direction for both registries.

\begin{table}[ht]
\centering
\footnotesize
\renewcommand{\arraystretch}{1.15}
\setlength{\tabcolsep}{3pt}
\resizebox{\textwidth}{!}{%
\begin{tabular}{l l cccc cccc}
\toprule
\multirow{2}{*}{Scenario} & \multirow{2}{*}{Model} & \multicolumn{4}{c}{Registry 1} & \multicolumn{4}{c}{Registry 2} \\
\cmidrule(lr){3-6} \cmidrule(lr){7-10}
 & & Tx (SD) & ES (SD) & $\Delta$Tx & $\Delta$ES & Tx (SD) & ES (SD) & $\Delta$Tx & $\Delta$ES \\
\midrule
\multirow{4}{*}{\makecell[l]{Tighter bound, Registry 2\\ $b_1=3$, $b_2=2$, $B=3$\\ $R_1,R_2\sim\mathrm{Pois}(5)$}}
 & Individual  & 50.42 (7.31) & 56.40 (2.62) & --      & --     & 46.72 (7.41) & 56.10 (2.65) & --      & --     \\
 & mKEP        & 56.52 (6.76) & 57.12 (2.47) & +12.1\% & +1.3\% & 52.48 (7.49) & 56.49 (2.63) & +12.3\% & +0.7\% \\
 & mKEP--IR    & 56.50 (7.56) & 56.86 (2.77) & +12.1\% & +0.8\% & 52.92 (7.73) & 56.07 (2.46) & +13.3\% & -0.1\% \\
 & mKEP--Fair  & 56.62 (7.39) & 56.73 (2.58) & +12.3\% & +0.6\% & 52.38 (8.12) & 55.88 (2.15) & +12.1\% & -0.4\% \\
\midrule
\multirow{4}{*}{\makecell[l]{Tighter bound, Registry 2\\ $b_1=3$, $b_2=2$, $B=3$\\ $R_1\sim\mathrm{Pois}(5)$, $R_2\sim\mathrm{Pois}(10)$}}
 & Individual  & 45.70 (9.24) & 56.37 (2.86) & --      & --     & 101.60 (11.34) & 56.36 (1.91) & --     & --     \\
 & mKEP        & 57.14 (8.58) & 57.39 (2.31) & +25.0\% & +1.8\% & 106.74 (12.02) & 56.04 (2.14) & +5.1\% & -0.6\% \\
 & mKEP--IR    & 55.68 (8.82) & 56.23 (2.27) & +21.8\% & -0.2\% & 107.68 (11.12) & 56.26 (1.75) & +6.0\% & -0.2\% \\
 & mKEP--Fair  & 54.56 (9.14) & 56.41 (2.24) & +19.4\% & +0.1\% & 107.52 (11.85) & 56.22 (1.68) & +5.8\% & -0.3\% \\
\midrule
\multirow{4}{*}{\makecell[l]{Tighter bound, Registry 2\\ $b_1=3$, $b_2=2$, $B=3$\\ $R_1\sim\mathrm{Pois}(10)$, $R_2\sim\mathrm{Pois}(5)$}}
 & Individual  & 107.54 (10.69) & 56.23 (1.52) & --     & --     & 44.92 (8.51) & 55.56 (2.59) & --      & --     \\
 & mKEP        & 111.46 (9.91)  & 55.80 (1.87) & +3.6\% & -0.7\% & 55.56 (7.10) & 55.77 (2.37) & +23.7\% & +0.4\% \\
 & mKEP--IR    & 113.36 (9.90)  & 55.92 (1.73) & +5.4\% & -0.5\% & 53.22 (7.66) & 55.89 (2.68) & +18.5\% & +0.6\% \\
 & mKEP--Fair  & 113.46 (9.70)  & 56.19 (1.88) & +5.5\% & -0.1\% & 52.60 (8.31) & 56.01 (2.51) & +17.1\% & +0.8\% \\
\bottomrule
\end{tabular}%
}
\caption{Effect of an asymmetric, tighter internal cycle-length bound for
Registry 2 ($b_1=3$ vs.\ $b_2=2$) together with a tightened cross-registry
bound ($B=3$), on transplant volume (Tx) and match quality (ES), crossed
against three arrival-rate combinations, by registry and allocation model.
Blood-group distribution is ASTRA for both registries throughout, isolating
the effect of the bound asymmetry from composition asymmetry. $\Delta$Tx and
$\Delta$ES report each pooled model's percentage change relative to that
registry's own Individual value in the same scenario, as in
Table~\ref{tab:results_arrival}. Standard deviations (in parentheses) are
computed across $N=50$ replications of 24 rounds each. Dropout probability
0.2 for both registries throughout.}
\label{tab:results_bound}
\end{table}

\subsection{Comparison with varying dropout probabilities}

The fourth factor influencing mKEP outcomes is the dropout probability governing each
participating registry - the likelihood that an unmatched pair exits the registry before the
following round. Dropout probabilities can vary due to several operational and logistical
factors, such as a patient's health deteriorating faster than other registry patients,
unwillingness to remain on the wait list due to longer wait times, donor refusal to donate, or
insufficient communication from the registry. Here, we compare the performance of mKEP to that of
individual registry solutions under different dropout-probability levels. The
analysis is based on blood group distributions from the ASTRA registry, with arrival rates held
at baseline (Poisson mean 5 for both registries) so that the effect of dropout probability can be
isolated from arrival-rate effects.

Table~\ref{tab:results_dropout} reports outcomes under a symmetric high-dropout configuration
($\mathrm{DP}_1=\mathrm{DP}_2=0.4$) and an asymmetric configuration
($\mathrm{DP}_1=0.2$, $\mathrm{DP}_2=0.4$), alongside the low-dropout symmetric baseline already
reported in Table~\ref{tab:results_arrival} ($\mathrm{DP}_1=\mathrm{DP}_2=0.2$). As dropout
probability rises from 0.2 to 0.4, standalone (Individual) transplant volume falls for both
registries, consistent with more pairs exiting the pool before they can be matched: Registry 1's
standalone volume falls from 49.46 at baseline to 42.86 under symmetric high dropout, and
Registry 2's from 49.32 to 47.02. The relative benefit of pooling rises correspondingly. At
baseline, gains under all three pooled mechanisms cluster tightly between 9.1\% and 10.3\% for
both registries. Under symmetric high dropout, gains widen to 10.5--18.2\%, and under dropout
asymmetry to 8.3--13.1\% -- both a wider range than baseline, and in the symmetric case
substantially higher. This is consistent with pooling directly counteracting the mechanism
driving the loss: a larger, combined pool offers more alternative matches per round, reducing the
number of rounds a pair must wait and hence its exposure to dropout, an effect that is naturally
larger when the standalone dropout risk is itself larger.

Within the asymmetric scenario, Registry 2 -- the registry with the higher dropout probability --
receives a larger gain than Registry 1 under both mKEP--IR (+11.3\% vs.\ +9.9\%) and mKEP--Fair
(+12.1\% vs.\ +8.3\%), consistent with the pattern observed under blood-group and cycle-bound
asymmetry elsewhere in this study, where the more structurally disadvantaged registry benefits
more from the fairness-oriented mechanisms. In the symmetric high-dropout scenario, both
mechanisms also give a larger gain to Registry 1 than to Registry 2 (+15.2\% vs.\ +11.7\% under
IR; +14.5\% vs.\ +12.7\% under Fair); however, since both registries share identical parameters
in this scenario, this split reflects the particular random realization underlying this set of
replications rather than a designed asymmetry, and should be interpreted with corresponding
caution. Match-quality changes are generally modest, with one exception: Registry 2's edge score
rises noticeably under dropout asymmetry (+3.2\% under mKEP, +2.5\% under mKEP--Fair), a larger
shift than observed for match quality in any other scenario examined in this study.

\begin{table}[ht]
\centering
\footnotesize
\renewcommand{\arraystretch}{1.15}
\setlength{\tabcolsep}{3pt}
\resizebox{\textwidth}{!}{%
\begin{tabular}{l l cccc cccc}
\toprule
\multirow{2}{*}{Scenario} & \multirow{2}{*}{Model} & \multicolumn{4}{c}{Registry 1} & \multicolumn{4}{c}{Registry 2} \\
\cmidrule(lr){3-6} \cmidrule(lr){7-10}
 & & Tx (SD) & ES (SD) & $\Delta$Tx & $\Delta$ES & Tx (SD) & ES (SD) & $\Delta$Tx & $\Delta$ES \\
\midrule
\multirow{4}{*}{\makecell[l]{High dropout, symmetric\\ $\mathrm{DP}_1=\mathrm{DP}_2=0.4$}}
 & Individual  & 42.86 (8.00) & 56.25 (2.72) & --      & --     & 47.02 (8.20) & 55.62 (2.64) & --      & --     \\
 & mKEP        & 50.66 (7.01) & 56.16 (3.28) & +18.2\% & -0.2\% & 51.96 (8.47) & 56.13 (2.53) & +10.5\% & +0.9\% \\
 & mKEP--IR    & 49.38 (7.85) & 56.45 (3.05) & +15.2\% & +0.4\% & 52.50 (8.21) & 56.32 (2.36) & +11.7\% & +1.3\% \\
 & mKEP--Fair  & 49.06 (8.27) & 55.89 (2.86) & +14.5\% & -0.6\% & 53.00 (8.81) & 55.90 (3.26) & +12.7\% & +0.5\% \\
\midrule
\multirow{4}{*}{\makecell[l]{Dropout asymmetry\\ $\mathrm{DP}_1=0.2$, $\mathrm{DP}_2=0.4$}}
 & Individual  & 50.92 (9.04) & 56.13 (2.81) & --      & --     & 43.62 (7.80) & 54.63 (3.11) & --      & --     \\
 & mKEP        & 57.60 (8.35) & 56.22 (2.12) & +13.1\% & +0.2\% & 47.80 (7.34) & 56.39 (2.71) & +9.6\%  & +3.2\% \\
 & mKEP--IR    & 55.94 (9.25) & 56.46 (2.28) & +9.9\%  & +0.6\% & 48.56 (7.58) & 55.07 (2.65) & +11.3\% & +0.8\% \\
 & mKEP--Fair  & 55.14 (9.04) & 56.38 (2.66) & +8.3\%  & +0.4\% & 48.90 (8.11) & 55.98 (2.88) & +12.1\% & +2.5\% \\
\bottomrule
\end{tabular}%
}
\caption{Effect of dropout-probability level and asymmetry on transplant volume
(Tx) and match quality (ES), by registry and allocation model. $\Delta$Tx and
$\Delta$ES report each pooled model's percentage change relative to that
registry's own Individual value in the same scenario, as in
Table~\ref{tab:results_arrival}. Standard deviations (in parentheses) are
computed across $N=50$ replications of 24 rounds each. Arrival rates and
blood-group distributions held at baseline ($R_1,R_2\sim\mathrm{Pois}(5)$,
ASTRA for both, $b_1=b_2=3$, $B=5$); the low-dropout symmetric case
($\mathrm{DP}_1=\mathrm{DP}_2=0.2$) is the Table~\ref{tab:results_arrival}
baseline and is not repeated here.}
\label{tab:results_dropout}
\end{table}

\begin{table}[ht]
\centering
\footnotesize
\renewcommand{\arraystretch}{1.15}
\setlength{\tabcolsep}{3pt}
\resizebox{\textwidth}{!}{%
\begin{tabular}{l l cccc cccc}
\toprule
\multirow{2}{*}{Scenario} & \multirow{2}{*}{Model} & \multicolumn{4}{c}{Registry 1} & \multicolumn{4}{c}{Registry 2} \\
\cmidrule(lr){3-6} \cmidrule(lr){7-10}
 & & Dropouts (SD) & Wait & $\Delta$Drop & $\Delta$Wait & Dropouts (SD) & Wait & $\Delta$Drop & $\Delta$Wait \\
\midrule
\multirow{4}{*}{\makecell[l]{High dropout, symmetric\\ $\mathrm{DP}_1=\mathrm{DP}_2=0.4$}}
 & Individual  & 70.08 (9.42) & 1.42 & --     & --     & 68.76 (10.48) & 1.41 & --     & --     \\
 & mKEP        & 62.96 (9.94) & 1.30 & -10.2\% & -8.7\% & 63.84 (9.02)  & 1.32 & -7.2\% & -6.7\% \\
 & mKEP--IR    & 64.22 (10.32) & 1.29 & -8.4\%  & -9.1\% & 63.22 (9.32)  & 1.28 & -8.1\% & -9.2\% \\
 & mKEP--Fair  & 64.00 (9.50)  & 1.28 & -8.7\%  & -9.6\% & 63.16 (9.61)  & 1.31 & -8.1\% & -7.4\% \\
\midrule
\multirow{4}{*}{\makecell[l]{Dropout asymmetry\\ $\mathrm{DP}_1=0.2$, $\mathrm{DP}_2=0.4$}}
 & Individual  & 60.56 (8.78) & 1.68 & --     & --      & 70.18 (9.00) & 1.45 & --     & --      \\
 & mKEP        & 55.12 (8.22) & 1.53 & -9.0\% & -8.7\%  & 66.40 (8.61) & 1.28 & -5.4\% & -11.5\% \\
 & mKEP--IR    & 56.24 (8.17) & 1.51 & -7.1\% & -9.9\%  & 65.84 (8.80) & 1.26 & -6.2\% & -13.0\% \\
 & mKEP--Fair  & 56.52 (8.01) & 1.48 & -6.7\% & -11.9\% & 65.58 (9.46) & 1.29 & -6.6\% & -10.8\% \\
\bottomrule
\end{tabular}%
}
\caption{Total dropouts over the 24-round horizon (mean, SD across $N=50$
replications in parentheses) and average waiting time in rounds until match,
for the same two dropout-probability scenarios and four allocation models as
Table~\ref{tab:results_dropout}. Waiting time is computed only over matched
patients; it does not include patients who dropped out unmatched. $\Delta$Drop
and $\Delta$Wait report each pooled model's percentage change relative to
Individual in the same scenario; as in Table~\ref{tab:results_arrival_dropwait},
\emph{negative} values here indicate improvement (fewer dropouts, shorter
waits) over acting alone.}
\label{tab:results_dropout_dropwait}
\end{table}

Comparing the dropout-probability scenarios in Table~\ref{tab:results_dropout_dropwait} against
the low-dropout baseline in Table~\ref{tab:results_arrival_dropwait} reveals a pattern that
initially appears counterintuitive: raising a registry's own dropout probability from 0.2 to 0.4
\emph{reduces} its average waiting time among matched patients, by roughly 12--17\%, and this
holds under every allocation model considered -- Individual, mKEP, mKEP--IR, and mKEP--Fair alike.
The same pattern appears in the asymmetric scenario: Registry 2's waiting time falls by
13--15\% when its own dropout probability is raised to 0.4, while Registry 1's waiting time is
essentially unchanged (0\% to $-2.6\%$) when only its partner's dropout probability rises and its
own remains at 0.2. This confirms that the effect is driven specifically by a registry's own
dropout risk rather than by any change in the pooled matching process itself.

This decrease in measured waiting time does not reflect patients being served faster; it reflects
how waiting time is defined. Because waiting time is computed only over patients who are
eventually matched, raising the dropout probability removes exactly the patients who would have
waited longest -- they exit the pool before a match can be found, and are therefore never counted
in the waiting-time average. The remaining, matched population is increasingly composed of
patients who happened to match quickly, which mechanically lowers the reported average even
though no patient is individually better off. Read together with the dropout counts in the same
table, which rise as expected when dropout probability increases, the two metrics describe two
sides of the same underlying effect: fewer patients survive long enough to be matched at all, and
the ones who do are, by construction, a faster-matching subset. Waiting time and dropout count
should therefore be interpreted jointly in this comparison; taken in isolation, the waiting-time
figures would misleadingly suggest that higher dropout probability improves outcomes, when the
opposite is true.

\section{Discussion and Conclusion} \label{sec: Conclusion}
 
This study set out to examine whether multi-registry kidney exchange (mKEP) reliably benefits
participating registries, and how that benefit is distributed when registries differ in arrival
rate, blood-group composition, cycle-length bound, and dropout probability. We compared four
allocation mechanisms -- standalone registry operation, unconstrained pooling, a cumulative
individual-rationality mechanism (mKEP--IR), and a Shapley-value-based fairness mechanism with a
rolling credit (mKEP--Fair) across twelve scenarios spanning these four sources of heterogeneity, using data and blood-group distributions calibrated to the Indian context.
 
Table~\ref{tab:results_benefit_ratio} collects the benefit ratio -- average transplants under
each pooled mechanism divided by that registry's own standalone average across all twelve scenarios. Four conclusions emerge from viewing the results this way, together.
 
Blood-group composition, not registry size or arrival rate, is what limits a registry's
incentive to join an mKEP. Under blood-group-composition asymmetry, the easier-to-match
registry's gain remains small even under both safeguarded mechanisms. It ranges as low as 0.6\%
under mKEP--IR and 1.8\% under mKEP--Fair when that registry also has the higher arrival rate
(Table~\ref{tab:results_bloodgroup}). So a fairness or individual-rationality criterion does
not, by itself, make participation clearly worthwhile for such a registry. The blood-type
breakdown in Table~\ref{tab:results_bloodgroup_bytype} traces this to O-type patients
specifically: the easier-to-match registry's O-type matches remain below their standalone level
in most blood-group-asymmetry scenarios even under mKEP--IR and mKEP--Fair, the one blood type
where the safeguards do not fully offset the composition-driven loss. Registries with an
easy-to-match donor pool should therefore expect blood-group composition, not size or
throughput, to be the main determinant of whether an mKEP is worth joining.
 
mKEP--Fair does not uniformly outperform mKEP--IR; which mechanism protects a registry better depends on what is constraining that registry. Under a tighter cycle-length bound specifically,
mKEP--IR gives a larger gain than mKEP--Fair to the more constrained registry in all three bound scenarios examined (Table~\ref{tab:results_bound}), a pattern that holds even when that registry additionally has the lower arrival rate. Under arrival-rate asymmetry alone, by contrast, the two mechanisms favor opposite sides: mKEP--IR gives the larger gain to the lower-arrival registry
while mKEP--Fair gives the larger gain to the higher-arrival one. Across all 24 registry-scenario combinations examined, mKEP--IR gives the larger gain in 14 cases and mKEP--Fair in 10 but one concentrated in the cycle-bound scenarios specifically rather than
distributed evenly across every source of heterogeneity. This is consistent with a structural
difference between the two mechanisms: mKEP--IR never revises a registry's credit downward once
earned, while mKEP--Fair claws back surplus in subsequent rounds to remain budget-balanced
(Section~\ref{sec: Model}). A registry designing an mKEP around a tighter internal cycle-length
bound may therefore be better served by a simple individual-rationality criterion, while a
registry operating under a more relaxed bound may find a fairness-based criterion the better
fit. Match quality, measured by edge score, is only modestly affected by either mechanism in almost every scenario examined (generally within about one percentage point of the Individual baseline).
 
Dropout probability is the one factor examined in this study where pooling's benefit is unambiguous. Unlike arrival-rate or blood-group asymmetry, where the more advantaged registry's
gain can shrink to near zero or turn negative, a higher dropout probability increases the
relative benefit of pooling for both registries under every mechanism
(Table~\ref{tab:results_dropout}): a larger, combined pool offers more alternative matches per
round, reducing the number of rounds a pair must wait and hence its exposure to dropout. This
makes dropout risk the clearest case in this study for participation in an mKEP being
straightforwardly advantageous, with no equity trade-off to weigh against it.
 
This pattern under the fairness-oriented mechanisms is not exclusively a blood-group phenomenon,
though blood-group asymmetry is where it occurs most often and most severely. Of the 48
registry-scenario-mechanism combinations under mKEP--IR and mKEP--Fair -- the two mechanisms a
registry would reasonably rely on, since both try to ensure it is not made worse off -- 7 yield a
gain below 5\% (italicized in Table~\ref{tab:results_benefit_ratio}), and all seven belong to
Registry 2. Five occur in the blood-group-asymmetry scenarios; the remaining two occur under
arrival-rate asymmetry alone, with no blood-group asymmetry present (Registry 2 at +3.2\% under
mKEP--IR and +4.6\% under mKEP--Fair when it has the higher arrival rate). A registry that is
both easy to match and already high-arrival can therefore find too little incentive to join an
mKEP for reasons beyond blood-group composition specifically.
\begin{table}[ht]
\centering
\footnotesize
\renewcommand{\arraystretch}{1.2}
\setlength{\tabcolsep}{4pt}
\resizebox{\textwidth}{!}{%
\begin{tabular}{l ccc ccc}
\toprule
 & \multicolumn{3}{c}{Registry 1 (Tx$_{\text{model}}$/Tx$_{\text{Ind}}$)} & \multicolumn{3}{c}{Registry 2 (Tx$_{\text{model}}$/Tx$_{\text{Ind}}$)} \\
\cmidrule(lr){2-4} \cmidrule(lr){5-7}
Scenario & mKEP & mKEP--IR & mKEP--Fair & mKEP & mKEP--IR & mKEP--Fair \\
\midrule
Baseline (5,5)                            & 1.103 & 1.091 & 1.092 & 1.095 & 1.103 & 1.097 \\
Arrival asymmetry (5,10)                  & 1.128 & 1.120 & 1.082 & \textit{1.035} & \textit{1.032} & \textit{1.046} \\
Arrival both high (10,10)                 & 1.065 & 1.059 & 1.054 & 1.052 & 1.057 & 1.051 \\
\addlinespace
BG asym., equal arrivals (5,5)            & 1.256 & 1.176 & 1.157 & \textbf{0.996} & \textit{1.045} & 1.052 \\
BG asym., unequal arrivals (5,10)         & 1.451 & 1.279 & 1.197 & \textbf{0.957} & \textit{1.006} & \textit{1.018} \\
BG asym., unequal arrivals (10,5)         & 1.137 & 1.087 & 1.092 & \textbf{0.991} & 1.062 & 1.052 \\
BG asym., equal arrivals (10,10)          & 1.244 & 1.135 & 1.103 & \textbf{0.945} & \textit{1.014} & \textit{1.026} \\
\addlinespace
Tighter bound ($b_2{=}2$), equal (5,5)    & 1.121 & 1.121 & 1.123 & 1.123 & 1.133 & 1.121 \\
Tighter bound ($b_2{=}2$), unequal (5,10) & 1.250 & 1.218 & 1.194 & 1.051 & 1.060 & 1.058 \\
Tighter bound ($b_2{=}2$), unequal (10,5) & \textit{1.036} & 1.054 & 1.055 & 1.237 & 1.185 & 1.171 \\
\addlinespace
High dropout, symmetric (0.4,0.4)         & 1.182 & 1.152 & 1.145 & 1.105 & 1.117 & 1.127 \\
Dropout asymmetry (0.2,0.4)               & 1.131 & 1.099 & 1.083 & 1.096 & 1.113 & 1.121 \\
\bottomrule
\end{tabular}%
}
\caption{Benefit ratio (average transplants under each pooled model, divided by
that registry's own Individual-model average) across every scenario examined
in this study. A ratio above 1 indicates the registry is better off pooling
under that mechanism than acting alone; a ratio above 1 but below 1.05 (\textit{italicized})
indicates a gain of less than 5\%, small enough that a registry may reasonably question whether
it justifies the coordination and data-sharing cost of participation; a ratio below 1
(\textbf{bolded}) indicates the registry is worse off pooling than acting alone. Comparing the
three ratios within a registry places how much of the unconstrained pool's
benefit survives under the two fairness-oriented mechanisms in direct view.
$R_1$ is ASTRA throughout; $R_2$ is APKD in the blood-group-asymmetry block
and ASTRA elsewhere. All values computed across $N=50$ replications of 24
rounds each.}
\label{tab:results_benefit_ratio}
\end{table}
 
This analysis has several limitations that should be acknowledged. First, the simulations relied
on blood-group distributions from only two registries -- ASTRA and APKD. While these provide
realistic and distinct settings, broader variation across registries may exist, and future work
could examine how mKEP outcomes change under additional or more diverse distributions. Second,
match quality in this study was assessed using only two parameters, HLA compatibility and
donor--recipient age difference, due to data constraints; incorporating additional clinical
indicators such as panel reactive antibody (PRA) levels, body mass index, or other immunological
factors could refine the assessment of match quality and potentially affect the observed
trade-offs between volume and quality. Third, the results reported here are based on a single
random seed per scenario; while standard deviations are reported throughout, formal statistical
testing of the differences between mechanisms -- particularly the smaller gains and the
mechanism-comparison findings above -- has not yet been performed, and would strengthen
confidence in the specific magnitudes reported.
 
Finally, this analysis does not incorporate non-directed donor (NDD) chains. This is a deliberate
scope decision rather than a data limitation: in India, non-directed donation and kidney exchange
are currently administered through separate institutional structures -- non-directed and deceased
donation are managed primarily through government-run programs, while kidney exchange programs of
the kind modeled in this study operate largely through private hospital networks. Modeling these
as separate systems, as we do here, reflects how mKEP actually operates in the Indian context
today, rather than a hypothetical integrated system. The formulation developed in this study can
be extended to incorporate NDD chains directly -- for instance, by introducing a non-directed
donor as a source node unconstrained by the cycle-length bound -- should India's institutional
structure evolve toward joint administration, or for researchers modeling settings where the two
mechanisms are already integrated. We leave this extension to future work.

\section*{Acknowledgments}

We gratefully acknowledge ASTRA for providing the data essential for this research. Their support
and cooperation were vital in advancing this study. Additionally, we extend our gratitude to
Dr. Peter Bir\'o for his suggestions as well.
 
\bibliographystyle{plainnat}
\bibliography{Reference}

\newpage

\setcounter{section}{0} 
\renewcommand{\thesection}{Appendix \Alph{section}} 

\end{document}